\documentclass[final]{siamart0516}
\usepackage[top=1.1in, bottom=1.3in, left=1.2in, right=1.2in]{geometry}
\usepackage{graphicx}
\usepackage{amsmath, amssymb, bm}
\usepackage{cases}
\usepackage{float, subfigure}
\usepackage{color}
\usepackage{amsfonts}
\usepackage{cite, url}
\usepackage[all,cmtip]{xy}
\usepackage{color,hyperref}
\definecolor{darkblue}{rgb}{0.0,0.0,0.6}
\hypersetup{colorlinks,breaklinks,
            linkcolor=darkblue,urlcolor=darkblue,
            anchorcolor=darkblue,citecolor=darkblue}

\newtheorem{assumption}{Assumption}
\newtheorem{remark}{Remark}

\def\0{{\bf 0}}
\def\1{{\bf 1}}

\def\beq{\begin{equation*}}
\def\eeq{\end{equation*}}
\def\bql{\begin{equation}}
\def\eql{\end{equation}}
\def\bqn{\begin{eqnarray*}}
\def\eqn{\end{eqnarray*}}
\def\bnl{\begin{eqnarray}}
\def\enl{\end{eqnarray}}
\def\bma{\begin{bmatrix}}
\def\ema{\end{bmatrix}}
\def\bmx{\begin{matrix}}
\def\emx{\end{matrix}}
\def\ben{\begin{enumerate}}
\def\een{\end{enumerate}}
\def\bit{\begin{itemize}}
\def\eit{\end{itemize}}
\def\bei{\begin{itemize}}
\def\eei{\end{itemize}}
\def\bet{\begin{tabular}}
\def\eet{\end{tabular}}

\newcommand{\ba}{\mathbf{a}}

\newcommand{\bb}{\mathbf{b}}
\newcommand{\bc}{\mathbf{c}}
\newcommand{\bd}{\mathbf{d}}
\newcommand{\bD}{\mathbf{D}}
\newcommand{\R}{\mathbb{R}}
\newcommand{\Gra}{\mathcal{G}^\mathrm{un}}
\newcommand{\Grad}{\mathcal{G}^\mathrm{dir}}
\newcommand{\GTIdir}{\mathcal{G}_\mathrm{TI}^\mathrm{dir}}
\newcommand{\GTIun}{\mathcal{G}_\mathrm{TI}^\mathrm{un}}
\newcommand{\GTVun}{\mathcal{G}_\mathrm{TV}^\mathrm{un}}
\newcommand{\GTVdir}{\mathcal{G}_\mathrm{TV}^\mathrm{dir}}
\newcommand{\V}{\mathcal{V}}
\newcommand{\E}{\mathcal{E}}
\newcommand{\A}{\mathcal{A}}
\newcommand{\bx}{\mathbf{x}}
\newcommand{\by}{\mathbf{y}}
\newcommand{\cbx}{\check{\bx}}
\newcommand{\cby}{\check{\by}}

\newcommand{\bq}{\mathbf{q}}
\newcommand{\br}{\mathbf{r}}
\newcommand{\bz}{\mathbf{z}}

\newcommand{\bs}{\mathbf{s}}
\newcommand{\f}{\mathbf{f}}

\newcommand{\barf}{f}
\newcommand{\barmu}{\bar{\mu}}
\newcommand{\hatmu}{\hat{\mu}}
\newcommand{\barL}{\bar{L}}
\newcommand{\bary}{\overline{y}}
\newcommand{\barx}{\overline{x}}
\newcommand{\baru}{\overline{u}}
\newcommand{\barkappa}{\bar{\kappa}}
\newcommand{\df}{\nabla\mathbf{f}}
\newcommand{\dfi}{\nabla f_i}

\newcommand{\dgi}{\nabla g_i}

\newcommand{\one}{\mathbf{1}}
\newcommand{\Fro}{\mathrm{F}}

\newcommand{\T}{\top}

\newcommand{\spa}[1]{\mathrm{span}\{#1\}}
\newcommand{\nul}[1]{\mathrm{null}\{#1\}}

\newcommand{\Wco}{\mathbf{C}}
\newcommand{\Wdo}{\mathbf{W}}

\newcommand{\wco}{C}
\newcommand{\wdo}{W}
\newcommand{\Ni}{\mathcal{N}_i}
\newcommand{\Niin}{\mathcal{N}_i^{\mathrm{in}}}
\newcommand{\Niout}{\mathcal{N}_i^{\mathrm{out}}}
\newcommand{\Njout}{\mathcal{N}_j^{\mathrm{out}}}
\newcommand{\Lap}{\textbf{{\L}}}

\newcommand{\sigmax}[1]{\sigma_{\max}\left\{#1\right\}}
\newcommand{\sig}{\delta}
\newcommand{\Cn}{J_1}
\newcommand{\Cnpush}{J_2}

\newcommand{\bu}{\mathbf{u}}
\newcommand{\bh}{\mathbf{h}}
\newcommand{\cbh}{\check{\bh}}
\newcommand{\bv}{\mathbf{v}}
\newcommand{\bV}{\mathbf{V}}

\newcommand{\tbR}{\widetilde{\mathbf{R}}}
\newcommand{\dia}[1]{\mathrm{diag}\left\{#1\right\}}
\newcommand{\BG}{B_{\ominus}}
\newcommand{\tBG}{\tilde{B}_{\ominus}}
\newcommand{\tB}{\tilde{B}}

\newcommand{\ttau}{\tilde{\tau}}
\newcommand{\tbd}{\tilde{\bd}}

\def\mr#1{{\color{black}#1}}

\newcommand{\TheTitle}{Achieving Geometric Convergence for Distributed Optimization over Time-Varying Graphs}
\newcommand{\TheAuthors}{Angelia Nedi\'c, Alex Olshevsky, and Wei Shi}

\headers{Geometric Convergence over Time-Varying Graphs}{\TheAuthors}

\title{{\TheTitle}\thanks{Initially submitted to the editors on July 13, 2016. Parts of the results have been submitted to the 55th IEEE Conference on Decision and Control and the 4th IEEE Global Conference on Signal and Information Processing for possible presentations. \funding{The work has been partially supported by NSF grant CNS 15-44953, Office of Naval Research under grant no. N00014-12-1-0998, and Air Force under grant number AF FA95501510394.}}}

\author{{\TheAuthors}\thanks{Dr. Angelia Nedi\'c and Dr. Wei Shi are with the School of Electrical, Computer, and Energy Engineering, Arizona State University, Tempe, AZ 85287, USA. Dr. Alex Olshevsky is with the Department of Electrical Engineering, Boston University, Boston, MA 02215, USA. (\email{Angelia.Nedich@asu.edu}, \email{Wilbur.Shi@asu.edu}, \email{alexols@bu.edu})}}

\begin{document}

\maketitle

\begin{abstract}
This paper considers the problem of distributed optimization over time-varying graphs. For the case of undirected graphs, we introduce a distributed algorithm, referred to as DIGing, based on a combination of a \underline{d}istributed \underline{i}nexact \underline{g}radient method and a gradient track\underline{ing} technique. The DIGing algorithm uses doubly stochastic mixing matrices and employs fixed step-sizes and, yet, drives all the agents' iterates to a global and consensual minimizer. When the graphs are directed, in which case the implementation of doubly stochastic mixing matrices is unrealistic, we construct an algorithm that incorporates the push-sum protocol into the DIGing structure, thus obtaining Push-DIGing algorithm. The Push-DIGing uses column stochastic matrices and fixed step-sizes, but it still converges to a global and consensual minimizer. Under the strong convexity assumption, we prove that the algorithms converge at R-linear (geometric) rates as long as the step-sizes do not exceed some upper bounds. We establish explicit estimates for the convergence rates. When the graph is undirected it shows that DIGing scales polynomially in the number of agents. We also provide some numerical experiments to demonstrate the efficacy of the proposed algorithms and to validate our theoretical findings.
\end{abstract}


\section{Introduction}\label{sec:intro}

This paper focuses on the following distributed convex optimization problem:
\begin{equation} \label{eq:F}
\min\limits_{x\in\R^p}~\barf(x)=\frac{1}{n}\sum\limits_{i=1}^n
f_i(x),
\end{equation}
where each function $f_i: \R^p\rightarrow \R$ is held privately by agent $i$ to encode the agent's objective function. We assume that the agents are connected through a communication network which can be time-varying. The agents want to collaboratively solve the problem, while each agent can only receive/send the information from/to its immediate neighbors (to be specified precisely soon). Problems of the form (\ref{eq:F}) that require distributed computing have appeared in various domains including information processing and decision making in sensor networks, networked vehicle/UAV coordination/control, as well as distributed estimation and learning. Some examples include distributed averaging \cite{Xiao2007,Cai2014,Olshevsky2014}, distributed spectrum sensing \cite{Bazerque2010}, formation control \cite{Ren2006,Olshevsky2010}, power system control \cite{Ram2009,Gan2013}, statistical inference and learning \cite{Rabbat2004,Forero2010,Nedic2015_2}. In general, distributed optimization framework fits the scenarios where the data is collected and/or stored in a network of agents and having a fusion center is either inapplicable or unaffordable. In such scenarios, data processing and computing is to be performed in a distributed but collaborative manner by the agents within the network.

We assume that the functions $f_i$ in problem \eqref{eq:F} are convex and continuously differentiable. For such a problem,  we propose a class of distributed algorithms that solve the problem over time-varying connectivity graphs for two different cases, namely, the case when the graphs are undirected and the case when they are directed. The algorithms employ consensus ideas for estimating the gradient of the global objective function in \eqref{eq:F}. When at least one of the objective functions is strongly convex, we show that the algorithms achieve R-linear convergence rates \footnote{Suppose that a sequence $\{x(k)\}$ converges to $x^*$ in some norm $\|\cdot\|$. We say that the convergence is: (i) Q-linear if there exists $\lambda\in(0,1)$ such that $\frac{\|x(k+1)-x^*\|}{\|x(k)-x^*\|} \leq\lambda$ for all $k$; (ii) R-linear if there exists $\lambda\in(0,1)$ and some positive constant $C$ such that $\|x(k)-x^*\|\leq C\lambda^k$ for all $k$. Both of these rates are geometric, and they are often referred to as global rates to be distinguished from the case when the given relations are valid for some sufficiently large indices $k$. The difference between these two types of geometric rate is in that Q-linear rate implies monotonic decrease of $\|x(k)-x^*\|$, while R-linear rate does not.}.

\subsection{Literature Review}
The research on distributed optimization dates back to 1980s \cite{Bertsekas1983,Tsitsiklis1986}. Since then, due to the emergence of large-scale networks, the development of distributed algorithms for solving problem in \eqref{eq:F} has received significant attention recently. \mr{Besides the decentralized and distributed approaches we are going to discuss below, many efforts have been made to solve \eqref{eq:F} in a master-slave structured isotropic network. The distributed algorithms designed over such special structure are usually fast in practice and mostly used in machine learning to handle big-data in a cluster of computers \cite{Boyd2011,Cevher2014}. Such scheme is ``centralized'' due to the use of a ``master''. In this paper, we focus on solving \eqref{eq:F} in a decentralized fashion motivated by the applications mentioned above.} 

Some earlier methods include distributed incremental (sub)gradient methods
\cite{Nedic2000,Nedic2001a,Nedic2001b,Ram2009_2} and incremental proximal methods \cite{Bertsekas2011,WangB2013}, while a more recent work includes incremental aggregated gradient methods~\cite{GOP2015} and its proximal gradient variants~\cite{Bertsekas2015}. All of the incremental methods require a special ring networks due to the nature of these methods. To handle a more general (possibly time-varying) networks, distributed subgradient algorithm was proposed in \cite{Nedic2009}, while its stochastic variant was studied in \cite{Ram2010} and its asynchronous variant in \cite{Nedic2011} with provable convergence rates.
These algorithms are intuitive and simple but usually slow due to the fact that even if the objective functions are differentiable and strongly convex, these methods still need to use diminishing step-size to converge to a consensual solution.
Other works on distributed algorithms that also require the use of diminishing step-sizes
include~\cite{Duchi2012,Zhu2012,Nedic2013}. With a fixed step-size, these distributed methods can be fast, but they only converge to a neighborhood of the solution set. This phenomenon creates an exactness-speed dilemma. A different class of distributed approaches that bypasses this dilemma is based on introducing Lagrangian dual variables and working with the Lagrangian function. The resulting algorithms include distributed dual decomposition~\cite{Terelius2011} and decentralized alternating direction method of multipliers (ADMM)~\cite{Bertsekas1997,Mateos2010}.
Specifically, the decentralized ADMM can employ a fixed step-size and it has nice provable rates~\cite{Wei2013}.
Under the strong convexity assumption, the decentralized ADMM has been shown to have linear convergence for time-invariant undirected graphs~\cite{Shi2014}. Building on (augmented) Lagrangian, a few improvements have been made via proximal-gradient \cite{Chang2014}, stochastic gradient \cite{Hong2015}, and second-order approximation \cite{Mokhtari2015,Mokhtari2016}. In particular, ADMM over a random undirected network has been shown to have $O(1/k)$ rate for convex functions~\cite{Hong2015}. However, de-synchronization and extensions of these methods to time-varying undirected graphs are more involved \cite{Wei2013,Hong2015}, while their extensions to directed graphs are non-existent in the current literature.

Some distributed methods exist that do not (explicitly) use dual variables but can still converge to an exact consensual solution while using fixed step-sizes. In particular, work in~\cite{Chen2012_2} employs multi-consensus inner loop and Nesterov's acceleration method, which gives a proximal-gradient algorithm with a rate at least $O(1/k)$. By utilizing multi-consensus inner loop, \emph{adapt-then-combine} (ATC) strategy, and Nesterov's acceleration, the algorithm proposed in \cite{Jakovetic2014} is shown to have $O\left(\ln(k)/k^2\right)$ rate under the assumption of bounded and Lipschitz gradients. For least squares, the general diffusion strategy (a generalization of ATC) can converge to the global minimizer`\cite{Sayed2013}. Although it is unknown to the literature, the above algorithms that do not use dual variable but use fixed step-size are not likely to reach linear convergence even under the strong convexity assumption. References \cite{Shi2015,Shi2015_2,Shi2015_3} use a difference structure to cancel the steady state error in decentralized gradient descent \cite{Nedic2009,Kun2014}, thereby developing the algorithm EXTRA and its proximal-gradient variant. EXTRA converges at an $o(1/k)$ rate when the objective function in \eqref{eq:F} is convex, and it has a Q-linear rate when the objective function is strongly convex.

Aside from the diminishing step-size issue, another topic of interest is distributed optimization over \emph{time-varying directed graphs}. The distributed algorithms over time-varying graphs require the use of doubly stochastic weight matrices, which are not easily constructed in a distributed fashion when the graphs are directed. To overcome this issue, reference \cite{Nedic2013} is the first to propose a different distributed approach, namely a subgradient-push algorithm that combines the distributed subgradient method \cite{Nedic2009} with the push-sum protocol \cite{Kempe2003}. While the subgradient-push eliminates the requirement of graph balancing \cite{Gharesifard2012}, it suffers from a slow sublinear\footnote{When an algorithm has convergence rate of $O(\theta(k))$, we say that the rate is sublinear if $\lim_{k\rightarrow+\infty}\frac{\lambda^k}{\theta(k)}=0$ for any constant $\lambda\in(0,1)$. A typical sublinear rates include $O(1/k^p)$ with $p>0$.} convergence rate even for strongly convex smooth functions due to its employment of diminishing step-size \cite{Nedic2014}. On the other hand, noticing that EXTRA has satisfactory convergence rates for undirected graphs, references \cite{Xi2015,Zeng2015} combine EXTRA with the push-sum protocol \cite{Kempe2003} to produce DEXTRA (ExtraPush) algorithm in hope of making it work over directed graph. It turns out that for a time-invariant strongly connected directed graph, DEXTRA converges at an R-linear rate under strong convexity assumption but the step-size has to be carefully chosen in some interval. However, the feasible set of step-sizes for DEXTRA can even be empty in some situations \cite{Xi2015}. \mr{Finally, the paper \cite{sun2016distributed} proposed an algorithm with diminishing step-size for  nonconvex optimization over directed graphs based on the push-sum method \cite{Kempe2003} and showed convergence to a stationary point.}

\mr{The algorithm we will study in this paper crucially relies on tracking differences of gradients and is a minor variation of algorithms that have appeared in \cite{Xu2015,Xu2016} and \cite{Lorenzo2016,Lorenzo2015,Lorenzo2016b}. To be specific, references \cite{Xu2015,Xu2016} utilize an Adapt-then-Combine variation \cite{Sayed2013} of the dynamic average consensus approach \cite{Zhu2010}, and thereby develop an Aug-DGM algorithm which is capable of employing uncoordinated step-sizes for multi-agent optimization. Independently, a scheme based on difference of gradients was proposed for the more general class of non-convex functions in~\cite{Lorenzo2016,Lorenzo2015,Lorenzo2016b} where a large class of distributed algorithms is developed.  Finally, we note that reference \cite{Qu2016}, appearing simultaneously with this work, also proposed a method for distributed optimization based on gradient differences. However, none of the papers mentioned in this paragraph provide a theoretical analysis for strongly convex optimization problems over {\it time-varying graphs}, which is the goal of this paper.}
\subsection{Summary of Contributions}

Prior to this work, an open question was how to construct a linearly convergent method for distributed optimization over time-varying (undirected or directed) graphs. 

\mr{The present paper resolves the issue. Specifically, we construct distributed methods which are linearly convergent over graphs which are time-varying and directed. Furthermore, we show that when the graphs are time-varying and {\em undirected}, a particular (distributed) choice of weights results in a polynomial iteration complexity (meaning that the number of iterations until the protocol reaches any fixed accuracy is polynomial in the total number of nodes).}

\mr{Our protocols work for all step-sizes small enough. To set the step-size to achieve the convergence rate guaranteed by our theorems, nodes need to know (i) an upper bound on the total number of agents and (ii) an upper bound on the number of time steps needed to achieve long-term connectivity. This compares favorably to some of the existing literature, e.g., \cite{Xi2015,Zeng2015}, 
which require detailed spectral information about the network for step-size selection.}

\mr{Moreover the technical tools we use are of independent interest. Although linearly convergent distributed optimization methods over {\em fixed} graphs were first developed in \cite{Shi2015}, extending the proof of \cite{Shi2015} to time-varying graphs does not appear to be possible. The current paper develops a new approach to the problem based on the small-gain theorem,  a standard tool for proving stability of interconnected dynamical systems in control theory. In fact, to our knowledge, our work is the first to use a small-gain based analysis to show convergence (and to bound convergence time) of an optimization protocol.}
\subsection{Paper Organization}

The rest of this paper is organized as follows. To facilitate the description of the technical ideas, the algorithms, and the analysis, we first introduce the notation in Subsection~\ref{sec:notation}. In Section~\ref{sec:undir_opt} we consider the case of undirected time-varying graphs, and we introduce a distributed consensus-based algorithm in Section~\ref{sec:algo_dev}. The algorithm uses ``\underline{d}istributed \underline{i}nexact \underline{g}radients'' and, also, employs a ``gradient track\underline{ing}'' technique, thus we term the algorithm as DIGing to account for its main design features. In Section \ref{sec:conv_analysis} we establish that the DIGing algorithm converges at an R-linear rate under standard assumptions including uniform joint strong connectivity of the graphs, the strong convexity of the objective function, and the Lipschitz continuity of the gradients. Moreover, we show that the convergence rate of DIGing scales polynomially in the total number of agents in the network. After this, we consider the case of time-varying directed graphs, and we propose a push-sum consensus-based variant of DIGing in Section~\ref{sec:dir_opt}. We establish its R-linear rate in Section \ref{sec:conv_analysis2}. Finally, some numerical simulations are given in~Section \ref{sec:num_exp}, and the paper concludes with some final remarks in~Section \ref{sec:concl}.

\subsection{Notation}\label{sec:notation}

Throughout the paper, the variable $x\in\R^p$ in the original problem~\eqref{eq:F} is viewed as a column vector. We let agent $i$ hold a \emph{local copy} of the variable $x$ of the problem in \eqref{eq:F}, which {is} denoted by $x_i\in\R^p$; its value at iteration/time $k$ is denoted by $x_i(k)$. We introduce an \emph{aggregate objective function} of the local variables: $\f(\bx)\triangleq\sum_{i=1}^{n} f_i(x_i)$, where its argument and gradient are defined as
\[
  \bx\triangleq\left(
     \begin{array}{ccc}
       \textrm{---}& x_1^\T & \textrm{---} \\
       \textrm{---}& x_2^\T & \textrm{---} \\
       &\vdots& \\
       \textrm{---}& x_n^\T & \textrm{---} \\
     \end{array}
   \right)\in\R^{n\times p}\quad\text{and}\quad
   \df(\bx)\triangleq\left(
     \begin{array}{ccc}
       \textrm{---}& (\nabla f_1(x_1))^\T & \textrm{---} \\
       \textrm{---}& (\nabla f_2(x_2))^\T & \textrm{---} \\
       &\vdots & \\
       \textrm{---}& (\nabla f_n(x_n))^\T & \textrm{---} \\
     \end{array}
   \right)\in\R^{n\times p},
\]
respectively. Each row $i$ of $\bx$ and $\df(\bx)$ is associated with agent $i$. We say that $\bx$ is \emph{consensual} if all of its rows are identical, i.e., $x_1=x_2=\cdots= x_n$. The analysis and results of this paper hold for all $p\geq1$. The reader can assume $p=1$ for convenience (so $\bx$ and $\df$ become vectors) without loss of generality. The notation is not standard but it enables us to present our algorithm and analysis in a compact form.

We let $\one$ denote a column vector with all entries equal to one (its size is to be understood from the context). For any matrix $\bv \in \R^{n\times p}$, we denote its average across the rows (a row corresponds to an agent) as $\overline{v} = \frac{1}{n}\bv^\T\one\in\R^{p}$, and its consensus violation as $\check{\bv} = \bv - \one{\overline{v}}^\T=\bv-\frac{1}{n}\one\one^\T\bv=\left(I-\frac{1}{n}\one\one^\T\right)\bv\triangleq\Lap\bv$, where $\Lap=I-\frac{1}{n}\one\one^\T$ is a symmetric matrix. For any $n\times p$ matrices $\ba$ and $\bb$, their inner product is denoted as $\langle\ba,\bb\rangle=\mathrm{Trace}(\ba^\T\bb)$. For a given matrix $\ba$, the Frobenius norm is given by $\|\ba\|_\Fro$, the (entry-wise) max norm is given by $\|\ba\|_{\max}$, while the spectral norm is given by $\|\ba\|_2$. $\|\ba\|_2$ equals to the largest singular value $\sigmax{\ba}$. We also use $\|\ba\|_{\Lap}$ to denote its $\Lap$ weighted (semi)-norm, that is, $\|\ba\|_{\Lap}=\sqrt{\langle\ba,\Lap\ba\rangle}$. Note that since $\Lap=\Lap^\T\Lap$, we always have $\|\ba\|_\Lap=\|\Lap\ba\|_\Fro.$ 
For any $n\times p$ matrices $\ba,\bb,\bc$ with $\bc=\ba+\bb$, in view of the definition of the consensus violation, it holds that
\[
\|\check\bc\|_\Fro = \|\bc\|_{\Lap} = \|\ba+\bb\|_{\Lap} \leq \|\ba\|_{\Lap}+\|\bb\|_{\Lap} = \|\check\ba\|_\Fro + \|\check\bb\|_\Fro,
\]
and we will use this property in the analysis without specific explanation. For any matrix $A\in\R^{m\times n}$, $\nul{A}\triangleq\{x\in\R^n\big|Ax=0\}$ is the null space of $A$ and $\spa{A}\triangleq \{y\in\R^m\big|y=Ax, x\in\R^n\}$ is the linear span of all the columns of $A$. 

\section{Distributed Optimization over Undirected Graphs}\label{sec:undir_opt}

In this section, we consider the case when the agents want to jointly solve the problem~\eqref{eq:F} while interacting over time-varying undirected graphs. We describe our algorithm, provide an interpretation of its steps and discuss its connection to some of the existing approaches. 

\subsection{The DIGing Algorithm}\label{sec:algo_dev}

We introduce the algorithm and provide some insights into its iterations. In what follows we will make use of the following proposition, which provides the optimality conditions for problem~\eqref{eq:F}.

\smallskip
\begin{proposition}[\cite{Shi2015}]\label{prop:conopt}
Assume $\nul{I-\Wdo}=\spa{\mathbf{1}}$ where $\Wdo\in\R^{n\times n}$. If $\bx^*$ satisfies conditions: \emph{(i)} $\bx^*=\Wdo\bx^*$ (consensus), \emph{(ii)} $\one^\T\df(\bx^*)=0$ (optimality), then the rows of $\bx^*$ are the same as each other and the transpose of each row is an optimal solution of the problem \eqref{eq:F}.
\end{proposition}
\smallskip

To illustrate the idea of the DIGing algorithm, let us focus on the case of a static graph for the moment. Consider the distributed gradient descent (DGD), given as follows:
\[
\label{eq:DGD_matrix}
\bx(k+1) = \Wdo\bx(k) - \alpha\df(\bx(k)),
\]
where $\Wdo$ is a doubly stochastic mixing matrix and $\alpha>0$ is a fixed step-size. The mixing part ``$\Wdo\bx(k)$'' is necessary for reaching consensus while DGD exhibits undesirable behavior due to its use of the gradient direction, ``$- \alpha\df(\bx(k))$''. To see this, let us break the update into steps per agent: for every agent $i$, we have $x_i(k+1)=\wdo_{ii}x_i(k)+\sum_{j\in\Ni}\wdo_{ij}x_j(k)-\alpha\dfi(x_i(k))$, where $\Ni$ is the set of the neighbors of agent $i$ in the given graph. Thus, each agent is updating using only the gradient of its local objective function $f_i$. Suppose now that the values $x_i(k)$ have reached consensus and that $x_i(k)=x^*$ for all $i$ and some solution $x^*$ of the problem~\eqref{eq:F}. Then, the mixing part gives $\wdo_{ii}x_i(k)+\sum_{j\in\Ni}\wdo_{ij}x_j(k)=x^*$ for all $i$. However, the gradient-based term gives $-\alpha\dfi(x_i(k))$ for all $i$, which need not be zero in general, thus resulting in $x_i(k+1)$ that will move away from the solution $x^*$ (recall that a solution to the problem \eqref{eq:F} is at a point $x$ where $\sum_{j=1}^n \nabla f_j(x)=0\ \forall i$ and not necessarily a point where $\dfi(x)=0\ \forall i$).

Conceptually, one (non-distributed) scheme that bypasses this limitation is the update
\begin{equation}\label{eq:DGD_ideal}
\bx(k+1) = \Wdo\bx(k) - \alpha \frac{1}{n}\one\one^\T \df(\bx(k))
\end{equation}
which can be implemented if every agent has access to the average of all the gradients $\nabla f_j(x_j(k))$, $j=1,\ldots,n$ (evaluated at each agent's local copy). One can verify that if \eqref{eq:DGD_ideal} converges, its limit point $\bx(\infty)$ satisfies the optimality conditions as given in Proposition \ref{prop:conopt}. However,
the update in \eqref{eq:DGD_ideal} is not distributed among the agents as it requires a central entity to provide the average of the gradients.

Nevertheless, one may approximate the update in \eqref{eq:DGD_ideal} through a surrogate direction that tracks the gradient average. To track the average of the gradients, namely, $\frac{1}{n}\one\one^\T\df(\bx(k))$, we introduce a variable $\by(k)$ that is updated as follows:
\begin{equation}\label{eq:DAC}
\by(k+1)=\Wdo\by(k)+\df(\bx(k+1))-\df(\bx(k)),
\end{equation}
with initialization $\by(0)=\df(\bx(0))$ and where each row $i$ of $\by(k)\in\R^{n\times p}$ is associated with agent $i$. A similar technique has been introduced in~\cite{Zhu2010} for dynamically tracking the average state of a multi-agent system, and for tracking some network-wide aggregate quantities in \cite{Ram2012,KNS2016}. If $\bx(k+1)$ converges to some point $\bx(\infty)$ and the underlying graph is connected, then it can be seen that the sequence $\by(k)$ generated by the gradient tracking procedure \eqref{eq:DAC} will converge to the point $\by(\infty)$ given by
\[
\by(\infty)=\frac{1}{n}\one\one^\T\df(\bx(\infty)),
\]
which is exactly what we need in view of \eqref{eq:DGD_ideal}. Replacing $\frac{1}{n}\one\one^\T\df(\bx(k))$ in \eqref{eq:DGD_ideal} by its dynamic approximation $\by(k)$ is exactly what we use to construct the DIGing algorithm. Furthermore, to accommodate time-varying graphs, the static weight matrix $\Wdo$ is replaced by a time varying matrix $\Wdo(k)$, thus resulting in the DIGing algorithm, as given below.

\smallskip
\begin{center}
  {\textbf{Algorithm 1: DIGing}}

  \smallskip
    \begin{tabular}{l}
    \hline
    \emph{  } Choose step-size $\alpha>0$ and pick any $\bx(0)\in \R^{n \times p}$;\\
    \emph{  } Initialize $\by(0)=\df(\bx(0))$;\\
    \emph{  } \textbf{for} $k=0,1,\ldots$ \textbf{do} \\
    \qquad$\bx(k+1)=\Wdo(k)\bx(k)-\alpha\by(k)$;\\
    \qquad$\by(k+1)=\Wdo(k)\by(k)+\df(\bx(k+1))-\df(\bx(k))$;\\
    \emph{  } \textbf{end}\\
    \hline
    \end{tabular}
\end{center}
\smallskip

Looking at an individual agent $i$, the initialization of DIGing uses an arbitrary $x_i(0)\in\R^p$ and sets $y_i(0)=\nabla f_i(x_i(0))$ for all $i=1,\ldots,n$. At each iteration $k$, every agent $i$ maintains two vectors, namely, $x_i(k),y_i(k)\in\R^p$, which are updated as follows:
\[
\begin{array}{l}
x_i(k+1)=\wdo_{ii}(k)x_i(k)+\sum_{j\in\Ni(k)}\wdo_{ij}(k)x_j(k)-\alpha y_i(k),\\
y_i(k+1)=\wdo_{ii}(k)y_i(k)+\sum_{j\in\Ni(k)}\wdo_{ij}(k)y_j(k)+\dfi(x_i(k+1))-\dfi(x_i(k)),
\end{array}
\]
where $\Ni(k)$ is the set of all neighbors of agent $i$ at time $k$. At every iteration $k$, each agent $i$ sends its current solution estimate $x_i(k)$ and average gradient estimate $y_i(k)$ to all of its neighbors $\Ni(k)$ while receiving all of its neighbors' solution estimates $x_j(k)$ and average gradient estimates $y_j(k)$, $\forall j\in\Ni(k)$. Then, each agent $i$ updates its vector $x_i(k+1)$ by mixing its own $x_i(k)$ and the neighbors' copies $x_j(k),$ for $j\in\Ni(k)$, with specific weights, and adjusting along the direction of $-y_i(k)$. Also, each agent $i$ updates its direction $y_i(k+1)$ by mixing its own $y_i(k)$ and the neighbors' directions $y_j(k)$, for $j\in\Ni(k)$ with specific weights, and by taking into account only the new information contained in the most recent gradient evaluation, as captured in the gradient-difference term $\dfi(x_i(k+1))-\dfi(x_i(k))$.

\subsection{Relation of DIGing to some of the existing approaches}

This section explains how the introduced DIGing algorithm is related to some of the other distributed algorithms.

\subsubsection{Connection with EXTRA \cite{Shi2015}}

When $\Wdo(k)=\Wdo$ (time-invariant case) and $\Wdo$ is symmetric, the DIGing algorithm  shares some similarity with EXTRA. If we eliminate the variables $\by(k)$ in the recursion of DIGing, we will obtain
\[
\bx(k+2)=(I+(2\Wdo-I))\bx(k+1)-\Wdo^2\bx(k)-\alpha[\df(\bx(k+1))-\df(\bx(k))].
\]
In this form, $(2\Wdo-I)$ and $\Wdo^2$ will be the two mixing matrices in EXTRA. As long as we have
\[
(2\Wdo-I)\preccurlyeq \Wdo^2\preccurlyeq\left(I+(2\Wdo-I)\right)/2\quad\text{and}\quad\Wdo^2\succ0,
\]
the convergence properties of DIGing will follow from the results in~\cite{Shi2015}. It can be seen that, when $\Wdo\succ0$, the convergence of DIGing follows from the convergence of EXTRA immediately. In this paper, we conduct the convergence analysis with more general choices of time-varying $\Wdo(k)$.

\subsubsection{Connection with primal-dual approaches}

DIGing has a primal-dual interpretation when the mixing matrices are static and symmetric. Indeed, suppose that $\Wdo(k)=\Wdo$ where $\Wdo$ is a symmetric doubly stochastic matrix, and consider the augmented Lagrangian function
\begin{equation}\label{eq:Lagrangian_GS}
\mathcal{L}_{\alpha}(\bx,\br)=\f(\bx)+\frac{1}{\alpha}\langle\br,(I-\Wdo)\bx\rangle+\frac{1}{2\alpha}\|\bx\|_{I-\Wdo^2}^2.
\end{equation}
If we apply the basic gradient method with a step-size $\alpha$ in Gauss-Seidel-like order for computing the saddle point of the augmented Lagrangian function \eqref{eq:Lagrangian_GS}, we will have
\[
\begin{array}{rcl}
\bx\text{-update: } \bx(k+1)&=&\bx(k)-\alpha\df(\bx(k))-(I-\Wdo)\br(k)-(I-\Wdo^2)\bx(k),\\
\br\text{-update: } \br(k+1)&=&\br(k)+(I-\Wdo)\bx(k+1).
\end{array}
\]
By eliminating the dual variable $\br$, we will obtain the same updates as in the DIGing algorithm 
(where the variable $\by$ is eliminated). The same will happen if, alternatively, we consider the augmented Lagrangian function
\begin{equation}\label{eq:Lagrangian_Jacobi}
\mathcal{L}_{\alpha}(\bx,\br)=\f(\bx)+\frac{1}{\alpha}\langle\br,(I-\Wdo)\bx\rangle+\frac{1}{2\alpha}\|\bx\|_{2I-2\Wdo}^2,
\end{equation}
and apply the basic gradient method with a step-size $\alpha$ in Jacobi-like order for seeking the saddle point of the augmented Lagrangian function \eqref{eq:Lagrangian_Jacobi}. \mr{Similar connections between the EXTRA algorithm and a general primal-dual have been made in a few recent references, \cite{Mokhtari2016}, \cite{mokhtari2016dsa}, and \cite{Hong2016}.}

\smallskip
\begin{remark}[\textbf{Primal-dual and symmetry}]
The above discussion assumes a time-invariant matrix $\Wdo$ which is symmetric. Even in the time-invariant graph case, but for asymmetric $\Wdo$, it appears to be difficult to adapt the classical primal-dual analysis for the recursion of DIGing. Our arguments in this paper are from a pure primal perspective and do not assume any symmetry property of $\Wdo(k)$. This suggests the possibility of its extension to the case of directed graphs, which will be addressed in Section \ref{sec:dir_opt}.
\end{remark}
\smallskip

\subsubsection{Connection with Aug-DGM \cite{Xu2015}}\label{sec:Xu2015}

A very recent reference that comes to our attention is \cite{Xu2015} that proposes a distributed consensus optimization algorithm, Aug-DGM, which is applicable to general time-invariant graphs and is quite similar to DIGing. The proposed algorithm is based on the combination of \emph{Adapt-then-Combine} (ATC) strategy of \cite{Sayed2013} and the dynamic average consensus for gradient tracking of\cite{Zhu2010}. It differs from DIGing only in the dynamic average gradient-consensus update, which uses an ATC variant. The updates of Aug-DGM are given by
\[
\begin{array}{rl}
\bx\text{-update: }& \bx(k+1)=\Wdo\left(\bx(k)-\bD\by(k)\right),\\
\by\text{-update: }& \by(k+1)=\Wdo\left(\by(k)+\df(\bx(k+1))-\df(\bx(k))\right),
\end{array}
\]
where $\Wdo$ is a doubly stochastic matrix and $\bD$ is a diagonal step-size matrix. When $\bD$ is chosen as $\alpha I$, it turns into an ATC variant of DIGing. With a general (positive) diagonal matrix $\bD$, Aug-DGM allows different agents to use different step-sizes and it still drives all the agent to reach a consensus on a global minimizer. The convergence of Aug-DGM is provided under general convexity and Lipschitz gradient assumptions.

\section{Convergence Analysis for DIGing over Undirected Graphs}\label{sec:conv_analysis}

In this section we establish the linear convergence of DIGing over time-varying undirected graphs. Let us formally describe the assumptions we make on the graphs and on the mixing matrices $\Wdo(k)$ that are compatible with the graphs.
Consider a time-varying undirected graph sequence $\{\Gra(0),\Gra(1),\ldots\}$. Every graph instance $\Gra(k)$ consists of a time-invariant set of agents $\V=\{1,2,\ldots,n\}$ and a set of time-varying edges $\E(k)$. The unordered pair of vertices ${(j,i)}\in\E(k)$ if and only if agents $j$ and $i$ can exchange information at time (iteration) $k$. The set of neighbors of agent $i$ at time $k$ is defined as $\Ni(k)=\left\{j\big|(j,i)\in\E(k)\right\}$.

In the sequel, we use the following notation:
\mr{\[\Gra_{b}(k)\triangleq\left\{\V,\E(k)\bigcup\E(k+1)\bigcup\cdots\bigcup\E(k+b-1)\right\}\ \text{ for any $k=0,1,\ldots$ and any $b=1,2,\ldots$},\]}
\[\Wdo_b(k)\triangleq\Wdo(k)\Wdo(k-1)\cdots\Wdo(k-b +1)\ \text{ for any $k=0,1,\ldots$ and any $b=0,1,\ldots$},\]
with the convention that $\Wdo_b(k)=I$ for any needed $k<0$ and $\Wdo_0(k)=I$ for any $k$.

Next, we give the basic assumption that we impose on the weight matrices.
\smallskip
\begin{assumption}[\textbf{Mixing matrix sequence $\{\Wdo(k)\}$}]\label{ass:matrix_W}
For any $k=0,1,\ldots$, the mixing matrix $\Wdo(k)=[\wdo_{ij}(k)]\in\R^{n\times n}$ satisfies the following relations:
\begin{itemize}
\item[\emph{(i)}] (Decentralized property) If $i\neq j$ and the edge $(j,i)\notin\E(k)$, then $\wdo_{ij}(k)=0$;
\item[\emph{(ii)}] (Double stochasticity) $\Wdo(k) \one = \one$, $\one^\T \Wdo(k) = \one^\T$;
\item[\emph{(iii)}] (Joint spectrum property) There exists a positive integer $B$ such that
\[
\sup_{k\ge B-1}\sig(k) < 1\ \text{ where }\ \sig(k)=\sigmax{\Wdo_B(k)-\frac{1}{n}\one\one^\T}\ \text{ for all }\ k=0,1,\ldots.
\]
\end{itemize}
\end{assumption}
\smallskip

In Assumption \ref{ass:matrix_W}, item (i) is due to the physical restriction of the network. Properties (ii) and (iii) are commonly used in the analysis of the rate of consensus algorithms. Several different mixing rules exist that yield the matrix sequences which have property (iii) (see subsection 2.4 of reference \cite{Shi2015}).

In particular, the following two assumptions taken together imply 
Assumption~\ref{ass:matrix_W}~\cite{Nedic2009_2}.

\smallskip
\mr{\begin{assumption}[\textbf{$\tB$-connected graph sequence}]\label{ass:connectedness}
The time-varying undirected graph sequence $\{\Gra(k)\}$ is $\tB$-connected. Specifically, there exists some positive integer $\tB$ such that the undirected graph $\Gra_{\tB}(t\tB)$ is connected for all $t=0,1,\ldots$.
\end{assumption}}
\smallskip

\mr{Assumption \ref{ass:connectedness} is typical for many results in multi-agent coordination and distributed optimization \cite{Nedic2015}. It is much weaker than the assumption of every $\Gra(k)$ being connected.}

\smallskip
\begin{assumption}[\textbf{Mixing matrix sequence $\{\Wdo(k)\}$}]\label{ass:matrix_W_specific}
For any $k=0,1,\ldots$, the mixing matrix $\Wdo(k)=[\wdo_{ij}(k)]\in\R^{n\times n}$ satisfies
\begin{itemize}
\item[\emph{(i)}] (Double stochasticity) $\Wdo(k) \one = \one$, $\one^\T \Wdo(k) = \one^\T$;
\item[\emph{(ii)}] (Positive diagonal) For all $i$, $\Wdo_{ii}(k) > 0$;
\item[\emph{(iii)}] (Edge utilization) If $(j,i) \in \E(k)$, then $\wdo_{ij}(k)>0$; otherwise $\wdo_{ij}(k)=0$;
\item[\emph{(iv)}] (Non-vanishing weights) There exists some $\tau>0$ such that if $\wdo_{ij}(k) > 0$, then $\wdo_{ij}(k) \ge\tau$;
\end{itemize}
\end{assumption}
\smallskip

Assumption \ref{ass:matrix_W_specific} is strong but typical for multi-agent coordination and optimization. For undirected graph it can be fulfilled, for example, by using Metropolis weights:
\[
\wdo_{ij}(k)=\left\{
         \begin{array}{ll}
           1/\left(1+\max\{d_i(k),d_j(k)\}\right),     &\text{ if } (j,i)\in\E(k),\\
           0,                              &\text{ if } (j,i)\notin\E(k)) \text{ and } j\neq i,\\
           1-\sum_{l\in\Ni(k)}\wdo_{il}(k),&\text{ if } j=i,
         \end{array}
       \right.
\]
where $d_i(k)=|\Ni(k)|$ is the degree of agent $i$ at time $k$. In this case, Assumption \ref{ass:matrix_W_specific} will be satisfied with the choice of $\tau=1/n$.

\mr{Now let us consider the graph sequence $\{\Gra(k)\}$ under Assumption \ref{ass:connectedness}. If the constant $B$ mentioned in Assumption \ref{ass:matrix_W} is chosen as $B\geq2\tB-1$, then for any $k=0,1,\ldots$, the union of a $B$-length consecutive clip of the edge set sequence from time index $k$, $\bigcup_{b=k}^{k+B-1}\E(b)$, is always a super set of $\bigcup_{b=\lceil k/\tB \rceil \tB}^{\lceil k/\tB \rceil \tB+\tB-1}\E(b)$ and the graph $\Gra_{\tB}(\lceil k/\tB \rceil \tB)$ by assumption is connected, thus it is very easy to see that the graph $\Gra_{B}(k)$ is connected. Here, the notation $\lceil\cdot\rceil$ is used to denote the ceiling function which rounds a real number to the least succeeding integer. With such a choice of $B$, Assumptions \ref{ass:connectedness} and \ref{ass:matrix_W_specific} together imply Assumption \ref{ass:matrix_W}. Regarding such a relationship, through out Section \ref{sec:conv_analysis}, we choose 
	\[B=2\tB-1\]
in the context wherever Assumption \ref{ass:connectedness} is used.}

The following lemma provides an important relation for later use.
\smallskip
\begin{lemma}[\textbf{$B$-step consensus contraction}]\label{lemma:mixing_contraction}
Under Assumption \ref{ass:matrix_W}, for any $k=B-1,B,\ldots$, and any matrix $\bb$ with appropriate dimensions, if $\ba = \Wdo_B(k) \bb$, then we have $\|\ba\|_{\emph{\Lap}} \leq \sig(k)\|\bb\|_{\emph{\Lap}}$, where $\sig(k)$ is as given in Assumption \ref{ass:matrix_W}(iii).
\end{lemma}
\smallskip

We do not claim any originality of this lemma. This lemma is fairly standard in consensus theory and it is a direct consequence of Assumption \ref{ass:matrix_W} due to the fact that $\Wdo_B(k)$ is doubly stochastic:
\[
\begin{array}{rcl}
\|\ba\|_{\Lap}
& =  &\|(I-\frac{1}{n}\one\one^\T)\Wdo_B(k)\bb\|_{\Fro}\\
& =  &\|(\Wdo_B(k)-\frac{1}{n}\one\one^\T)(I-\frac{1}{n}\one\one^\T)\bb\|_\Fro\\
&\leq&\sigmax{\Wdo_B(k)-\frac{1}{n}\one\one^\T}\|\bb\|_\Lap.
\end{array}
\]
To make our arguments more concise, we will use $\sig=\sup_{k\ge B-1}\{\sig(k)\}$ in our analysis of the algorithm. An explicit expression of $\delta$ in terms of $n$ can be found in \cite{Nedic2009_2} if the more specific Assumption \ref{ass:matrix_W_specific} is made.

We also need the following two assumptions on the objective functions, which are standard for deriving linear (geometric) rate of gradient algorithms for minimizing strongly convex smooth functions.
\smallskip
\begin{assumption}[\textbf{Smoothness}]\label{ass:smooth} For every agent $i$, its objective $f_i:\R^p\rightarrow\R$ is differentiable and has Lipschitz continuous gradients, i.e., 
there exists a Lipschitz constant $L_i\in(0,+\infty)$ such that
\[
\|\nabla f_i(x) - \nabla f_i(y)\|_\Fro \leq L_i \|x - y\|_\Fro\ \text{ for any }\ x, y \in \R^p.
\]
\end{assumption}
\smallskip

When Assumption~\ref{ass:smooth} holds, we will also say that each $\nabla f_i$ is $L_i$-Lipschitz (continuous). In the forthcoming analysis, we will use $L\triangleq\max_i \{L_i\}$, which is the Lipschitz constant of $\df(\bx)$, and $\barL\triangleq (1/n)\sum_{i=1}^n L_i$ which is the Lipschitz constant of $\nabla f(x)$.

\smallskip
\begin{assumption}[\textbf{Strong convexity}]\label{ass:strongly_convex}
For every agent $i$, its objective $f_i:\R^p\rightarrow\R$ satisfies
\[
f_i(x)\geq f_i(y)+\langle\dfi(y),x -y\rangle+\frac{\mu_i}{2}\|x-y\|_\Fro^2\ \text{ for any }\ x, y \in \R^p,
\]
where $\mu_i\in[0,+\infty)$ and at least one $\mu_i$ is nonzero.
\end{assumption}
\smallskip

When $\mu_i>0$, 
we will say that $f_i$ is $\mu_i$-strongly convex. In the analysis we will use $\hatmu\triangleq\max_i\{\mu_i\}$ and $\barmu\triangleq(1/n)\sum_{i=1}^n \mu_i$. Assumption \ref{ass:strongly_convex} implies the $\barmu$-strong convexity of $f(x)$. Under this assumption, the optimal solution to problem \eqref{eq:F} is guaranteed to exist and to be unique since $\barmu>0$. We note that all the convergence results in our analysis are achieved under Assumption \ref{ass:strongly_convex}. We will also use $\barkappa\triangleq L/\barmu$.

To establish the R-linear rate of the algorithm, one of our technical innovations will be to resort to a somewhat unusual version of small gain theorem under a well-chosen metric, whose original version has received an extensive research and been widely applied in control theory \cite{Desoer2009}. We will give an intuition of the whole analytical approach shortly, after stating the small gain theorem at first.

\subsection{The Small Gain Theorem}

Let us adopt the notation $\bs^i$ for the infinite sequence $\bs^i=\left(\bs^i(0),\bs^i(1),\bs^i(2),\ldots\right)$ where $\bs^i(k)\in\R^{n\times p},\ \forall i$. Furthermore, let us define
\begin{equation}\label{eq:norm}
\|\bs^i\|_\Fro^{\lambda, K} \triangleq \max_{k=0,\ldots,K} \frac{1}{\lambda^k}  \|\bs^i(k)\|_\Fro\quad\text{and}\quad\|\bs^i\|_\Fro^{\lambda} \triangleq \sup_{k \geq 0} \frac{1}{\lambda^k} \|\bs^i(k)\|_\Fro,
\end{equation}
where the parameter $\lambda\in(0,1)$ will serve as the linear rate parameter later in our analysis. While $\|\bs^i\|_\Fro^{\lambda, K}$ is always finite, $\|\bs^i\|_\Fro^{\lambda}$ may be infinite. If $n=p=1$, i.e., each $\bs^i(k)$ is a scalar, we will just write $|\bs^i|^{\lambda, K}$ and $|\bs^i|^{\lambda}$ for these quantities. Intuitively, $\|\bs^i\|_\Fro^{\lambda, K}$ is a weighted ``ergodic norm'' of $\bs^i$. Noticing that the weight $\frac{1}{\lambda^k}$ is exponentially growing with respect to $k$, if we can show that $\|\bs^i\|_\Fro^{\lambda}$ is bounded, then it
would imply that $\|\bs^i(k)\|_\Fro\to0$ geometrically fast. This ergodic definition enables us to give analysis to those algorithms which do not converge Q-linearly. Next we will state the small gain theorem which gives a sufficient condition to for the boundedness of $\|\bs^i\|_\Fro^{\lambda}$. The theorem is a basic result in control systems and a detailed discussion about its result can be found in \cite{Desoer2009}. For the sake of completeness, we include the proof.

\begin{theorem}[\textbf{The small gain theorem}]\label{theorem:small} Suppose that $\bs^1, \ldots, \bs^m$ are sequences such that for all positive integers $K$ and for each $i=1, \ldots, m$,  we have an arrow $\bs^i\rightarrow\bs^{(i \mod m)+1}$, that is
\begin{equation}\label{eq:smallgainrecur}
\|\bs^{(i~{\rm mod}~m)+1} \|_\Fro^{\lambda, K} \leq \gamma_i  \|\bs^{i} \|_\Fro^{\lambda, K} + \omega_i,
\end{equation}
where the constants (gains) $\gamma_1, \ldots, \gamma_m$ are nonnegative and satisfy $\gamma_1 \gamma_2 \cdots \gamma_m < 1$. Then
\[
\begin{array}{rcl}
\|\bs^1\|_\Fro^{\lambda}
&\leq&\frac{1}{1 - \gamma_1\gamma_2 \cdots \gamma_m}\left(\omega_1\gamma_2\gamma_3 \cdots \gamma_m + \omega_2\gamma_3\gamma_4 \cdots \gamma_m + \cdots + \omega_{m-1} \gamma_m  +  \omega_m\right).
\end{array}
\]
\end{theorem}
\smallskip
\begin{proof}
By iterating inequality \eqref{eq:smallgainrecur} for $i$ from $m$ down to $1$, we obtain
\[
\begin{array}{rcl}
\|\bs^1\|_\Fro^{\lambda, K}
&\leq& \gamma_m\gamma_{m-1} \cdots \gamma_1  \|\bs^1\|_\Fro^{\lambda, K} + \gamma_m\gamma_{m-1}\cdots\gamma_2\omega_1\\
&    & + \gamma_m\gamma_{m-1}\cdots\gamma_3\omega_2 + \cdots + \gamma_m\omega_{m-1} + \omega_m.
\end{array}
\]
Thus,
\begin{equation}\label{eq:sg_proof1}
\begin{array}{rcl}
\|\bs^1\|_\Fro^{\lambda, K}
&\leq & \frac{1}{1 - \gamma_1\gamma_2 \cdots \gamma_m} \left( \omega_1 \gamma_2\gamma_3 \cdots \gamma_m + \omega_2 \gamma_3\gamma_4 \cdots \gamma_m + \cdots + \omega_{m-1} \gamma_m  +  \omega_m\right).
\end{array}
\end{equation}
Since \eqref{eq:sg_proof1} holds for all $K$ and its right-hand side does not depend on $K$, taking $K\rightarrow\infty$ implies the desired relation.
\end{proof}
\smallskip

Clearly, the small gain theorem involves a cycle $\bs^1\rightarrow\bs^2\rightarrow\cdots\rightarrow\bs^m\rightarrow\bs^1$. Due to this cyclic structure [cf. \eqref{eq:smallgainrecur}], similar bounds hold for $\|\bs^i\|_\Fro^{\lambda},\ \forall i$.

\smallskip
\begin{lemma}[\textbf{Bounded norm $\Rightarrow$ R-linear rate}]\label{lemma:bound_R-linear}
For any matrix sequence $\bs^i$, if $\|\bs^i\|_\Fro^{\lambda}$ is bounded, then $\|\bs^i(k)\|_\Fro$ converges at a global R-linear (geometric) rate $O(\lambda^k)$.
\end{lemma}
\smallskip
\begin{proof}
If $\|\bs^i\|_\Fro^{\lambda}\leq U$ where $U$ is some nonnegative constant, then by the definition we have $\sup_{k \geq 0}~ \frac{1}{\lambda^k} \|\bs^i(k)\|_\Fro\leq U$, thus $\frac{1}{\lambda^k} \|\bs^i(k)\|_\Fro\leq U,\ \forall k$. The conclusion follows immediately from $\|\bs^i(k)\|_\Fro\leq U\lambda^k,\ \forall k$.
\end{proof}
\smallskip

\subsection{Sketch of the Main Idea}

Before summarizing our main proof idea, let us define some quantities which we will use frequently in our analysis. We define $\bx^*\triangleq\one(x^*)^\T$ where $x^*$ is the optimal solution of problem \eqref{eq:F}. Also, define
\[
\bq(k)\triangleq\bx(k) - \bx^*\ \text{ for any }\ k=0,1,\ldots,
\]
which  is the optimality residual of the iterates $\bx(k)$ (at the $k$-th iteration). Moreover, let us adopt the notation
\[
\bz(k)\triangleq\df(\bx(k)) - \df(\bx(k-1))\ \text{ for any }\ k=1,2,\ldots,
\]
and with the convention that $\bz(0)\triangleq\0$.

We will apply the small gain theorem with the $\|\cdot\|_\Fro^{\lambda, K}$ metric and a right choice of $\lambda < 1$ around the following circle of arrows:
\begin{equation}\label{eq:cycle_alg1}
\text{Algorithm 1: }\bq \rightarrow \bz \rightarrow \cby \rightarrow \cbx \rightarrow \bq,
\end{equation}
where, recall, $\bq$ is the difference between local copies and the global optimizer, $\bz$ is the successive difference of gradients, $\cby$ is the consensus violation of the estimation of gradient average across agents, and $\cbx$ is the consensus violation of local copies (see Subsection \ref{sec:notation} for the definition of operator ``$\check{\ }$'').

Intuitively: as long $\bq$ is small, the successive difference of the gradients $\bz$ is small since the gradients are close to zero in the neighborhood of the optimal point; as long as the successive difference of the gradients $\bz$ is small, the structure of DIGing implies that $\by$ is close to consensual; as long as $\by$ is close to consensual, then by the structure of DIGing so is $\bx$; and, finally, as long as $\bx$ is close to consensual, DIGing is very similar to gradient descent and drives the distance to the optimal point $\bq$ to zero and thus completes the cycle.

After the establishment of each arrow, we will apply the small gain theorem to conclude that every corresponding quantities under the metric $\|\cdot\|_\Fro^{\lambda, K}$ is bounded and hence conclude that all quantities in the ``circle of arrows'' decay at an R-linear rate $O(\lambda^k)$.

Note that to apply the small gain theorem, we would need to have gains ($\gamma_i$) that multiply to less than one. This is achieved by choosing an appropriate step-size $\alpha$. Indeed, by looking at the algorithm, we can see that the step-size appears only in one place -- the third arrow (i.e., the arrow $\cby \rightarrow \cbx$), and the dependence of the corresponding gain in that arrow is \emph{linear} in $\alpha$. Thus we should be able to apply the small gain theorem after choosing small enough $\alpha$.

\subsection{The Establishment of Each Arrow}

We now discuss the establishment of each arrow/relation in the sketch above [cf. \eqref{eq:cycle_alg1}].

The first arrow demonstrated in Lemma \ref{lemma:first_arrow} is a simple consequence of Assumption \ref{ass:smooth} (namely, it is a consequence of the fact that the gradient of ${\bf f}$ is $L$-Lipschitz).

\smallskip
\begin{lemma}[\textbf{Algorithm 1: The first arrow $\bq\rightarrow \bz$}]\label{lemma:first_arrow}
Under Assumption \ref{ass:smooth}, we have that for all $K=0,1,\ldots$ and any $\lambda\in (0,1)$,
\[
\|\bz\|_\Fro^{\lambda, K} \leq L \left( 1 + \frac{1}{\lambda} \right) \|\bq\|_\Fro^{\lambda, K}.
\]
\end{lemma}
\smallskip
\begin{proof}
By Assumption \ref{ass:smooth}, $\df(\bx)$ is $L$-Lipschitz and we have
\begin{equation}\label{eq:first_arrow_proof1}
\begin{array}{rcl}
\|\df(\bx(k+1)) - \df(\bx(k))\|_\Fro
& \leq & L \|\bx(k+1) - \bx(k)\|_\Fro \\
&  =   & L \|( \bx(k+1) - \bx^*) - ( \bx(k) - \bx^* ) \|_\Fro \\
& \leq & L \|\bx(k+1) - \bx^*\|_\Fro + L \|\bx(k) - \bx^*\|_\Fro.
\end{array}
\end{equation}
By the definition of $\bz$ and $\bq$, it follows from \eqref{eq:first_arrow_proof1} that
\begin{equation}\label{eq:first_arrow_proof2}
\begin{array}{rcl}
\lambda^{-(k+1)} \|\bz(k+1)\|_\Fro
&\leq& L \lambda^{-(k+1)} \| \bq(k+1)\|_\Fro + \frac{L}{\lambda} \lambda^{-k} \|\bq(k)\|_\Fro.
\end{array}
\end{equation}
Taking $\max_{k=0,1,\ldots,K-1}\{\cdot\}$ on both sides of \eqref{eq:first_arrow_proof2} gives
\[
\begin{array}{l}
\|\bz\|_\Fro^{\lambda, K}
 \leq  L \|\bq\|_\Fro^{\lambda, K}+ \frac{L}{\lambda} \|\bq\|_\Fro^{\lambda, K-1}
 \leq  L\left(1+\frac{1}{\lambda}\right)\|\bq\|_\Fro^{\lambda, K}.
\end{array}
\]
\end{proof}
\smallskip

Next we provide the lemmata for the second and third arrows in the cycle \eqref{eq:cycle_alg1}. They are proved by an almost identical analysis based on Lemma \ref{lemma:mixing_contraction}: indeed, a glance at the structure of DIGing implies that some (semi)-norm of $\bx$ can be bounded in terms of some (semi)-norm of $\by$, while some (semi)-norm of $\by$ can be bounded in terms of some (semi)-norm of $\bz$. This is a fairly straightforward application of Lemma \ref{lemma:mixing_contraction}, which shows how multiplication by $\Wdo(k)$ shrinks the distance toward the consensus subspace.

\smallskip
\begin{lemma}[\textbf{Algorithm 1: The second arrow $\bz\rightarrow\cby$}]\label{lemma:second_arrow}
Let Assumption \ref{ass:matrix_W} hold, and let $\sig=\sup_{k\ge B-1}\{\sig(k)\}$, where $\sig(k)$ is as given in Assumption \ref{ass:matrix_W}(iii). Also, let $\lambda$ be such that $\sig < \lambda^B < 1$. Then, we have for all $K=0,1,\ldots$,
\begin{equation}\label{eq:second_arrow_alg1}
\|\cby\|_\Fro^{\lambda,K}\leq\frac{\lambda(1-\lambda^B)}{(\lambda^B-\sig)(1-\lambda)}\|\bz\|_\Fro^{\lambda,K}
+\frac{\lambda^B}{\lambda^B-\sig}\sum\limits_{t=1}^B\lambda^{1-t}\|\cby(t-1)\|_{\Fro}.
\end{equation}
\end{lemma}
\smallskip
\begin{proof}
The equivalent relation in DIGing involving $\by$ and $\bz$ is
\begin{equation}\label{eq:second_arrow_proof0}
\by(k+1)=\Wdo(k)\by(k)+\bz(k+1).
\end{equation}
From \eqref{eq:second_arrow_proof0}, using Lemma \ref{lemma:mixing_contraction}, for all $k \geq B-1$, it follows that
\begin{equation}
\begin{array}{rcl}\label{eq:second_arrow_proof1_prime}
\|\cby(k+1)\|_\Fro
&  =   &\|\by(k+1)\|_\Lap\\
& \leq &\left\|\Wdo_B(k)\by(k+1-B)\right\|_{\Lap} + \left\|\Wdo_{B-1}(k)\bz(k+2-B)\right\|_{\Lap}\\
&      &+ \cdots + \|\Wdo_1(k)\bz(k)\|_{\Lap} + \|\Wdo_0(k)\bz(k+1)\|_{\Lap}\\
& \leq &\sig\|\cby(k+1-B)\|_\Fro + \sum\limits_{t=1}^B\|\bz(k+2-t)\|_\Fro,
\end{array}
\end{equation}
and therefore, for all $k=B-1,B,\ldots$,
\begin{equation}\label{eq:second_arrow_proof1}
\hspace{-1em}
\begin{array}{c}
\lambda^{-(k+1)}\| \cby(k+1)\|_\Fro
\leq\frac{\sig}{\lambda^{B}}\lambda^{-(k+1-B)} \|\cby(k+1-B)\|_\Fro+\sum\limits_{t=1}^B\frac{1}{\lambda^{t-1}}\lambda^{-(k+2-t)}\|\bz(k+2-t)\|_\Fro.
\end{array}
\end{equation}
To utilize the norm $\|\cdot\|_\Fro^{\lambda,K}$, we need to take $\max_{k=0,...,K}$, which in turn requires a relation for $\lambda^{-(k+1)}\| \cby(k+1)\|_\Fro$ with $k<B-1$. To obtain such a relation, we complement  the initial relation for \eqref{eq:second_arrow_proof1}, i.e.,
\begin{equation}\label{eq:second_arrow_proof1_init}
\begin{array}{c}
\lambda^{-(k+1)}\| \cby(k+1)\|_\Fro
\leq\lambda^{-(k+1)}\| \cby(k+1)\|_\Fro
\end{array}
\end{equation}
for $k=-1,\ldots,B-2$. Taking the maximum over $k=-1,0,\ldots,B-2$ on both sides of \eqref{eq:second_arrow_proof1_init} and the maximum over $k=B-1,\ldots, K$ in \eqref{eq:second_arrow_proof1}, and then by combining the obtained relations, we obtain
\[\label{eq:second_arrow_proof2}
\begin{array}{rcl}
\|\cby\|_\Fro^{\lambda,K}
&\leq&\frac{\sig}{\lambda^B}\|\cby\|_\Fro^{\lambda,K-B}
+\sum\limits_{t=1}^B\frac{1}{\lambda^{t-1}}\|\bz\|_\Fro^{\lambda,K+1-t}
+\sum\limits_{t=1}^B\lambda^{1-t}\|\cby(t-1)\|_{\Fro}\\
&\leq&\frac{\sig}{\lambda^B}\|\cby\|_\Fro^{\lambda,K}
+\sum\limits_{t=1}^B\frac{1}{\lambda^{t-1}}\|\bz\|_\Fro^{\lambda,K}
+\sum\limits_{t=1}^B\lambda^{1-t}\|\cby(t-1)\|_{\Fro}.
\end{array}
\]
Hence,
\[
\begin{array}{rcl}
\|\cby\|_\Fro^{\lambda,K}
&\leq&\frac{\lambda^B\sum\limits_{t=1}^B\frac{1}{\lambda^{t-1}}}{\lambda^B-\sig}\|\bz\|_\Fro^{\lambda,K}+\frac{\lambda^B}{\lambda^B-\sig}\sum\limits_{t=1}^B\lambda^{1-t}\|\cby(t-1)\|_{\Fro}\\
& =  &\frac{\lambda(1-\lambda^B)}{(\lambda^B-\sig)(1-\lambda)}\|\bz\|_\Fro^{\lambda,K}+\frac{\lambda^B}{\lambda^B-\sig}\sum\limits_{t=1}^B\lambda^{1-t}\|\cby(t-1)\|_{\Fro},
\end{array}
\]
which is exactly \eqref{eq:second_arrow_alg1}.
\end{proof}
\smallskip

\begin{lemma}[\textbf{Algorithm 1: The third arrow $\cby\rightarrow\cbx$}]\label{lemma:third_arrow}
Let Assumption \ref{ass:matrix_W} hold, and let $\sig=\sup_{k\ge B-1}\{\sig(k)\}$, where $\sig(k)$ is as given in Assumption \ref{ass:matrix_W}(iii). Furthermore, let $\lambda$ be such that $\sig < \lambda^B < 1$. Then, we have for all $K=0,1,\ldots$,
\[
\|\cbx\|_\Fro^{\lambda, K}\leq\frac{\alpha(1-\lambda^B)}{(\lambda^B-\sig)(1-\lambda)} \|\cby\|_\Fro^{\lambda, K}+ \frac{\lambda^B}{\lambda^B-\sig}\sum\limits_{t=1}^{B}\lambda^{1-t}\|\cbx(t-1)\|_\Fro.\label{eq:third_arrow_alg1}
\]
\end{lemma}
\smallskip

The relation in DIGing involving $\bx$ and $\by$ is given by
\begin{equation}\label{eq:third_arrow_proof0}
\bx(k+1)=\Wdo(k)\bx(k)-\alpha\by(k).
\end{equation}
Noticing the similarity between \eqref{eq:third_arrow_proof0} and \eqref{eq:second_arrow_proof0}, we omit the proof of Lemma \ref{lemma:third_arrow} since it is almost identical to that of Lemma \ref{lemma:second_arrow}.

With all the above lemmata in place, the last arrow of our proof sketch remains to be addressed. For this, we need an interlude on gradient descent with errors in the gradient. Since this part is relatively independent from the preceding development, we provide it in the next subsection.

\subsection{The Inexact Gradient Descent on a Sum of Strongly Convex Functions}

In this subsection, we consider the basic (centralized) first-order method for problem \eqref{eq:F} under inexact first-order oracle.  To distinguish from the notation used for our distributed optimization problem/algorithm/analysis, let us make some definitions that are only used in this subsection. Problem \eqref{eq:F} is restated as follows with different notation,
\[\label{eq:F_diff}
\min\limits_{x\in\R^d} g(x)=\frac{1}{n}\sum\limits_{i=1}^n g_i(x),
\]
where all $g_i$'s satisfy Assumptions \ref{ass:smooth} and \ref{ass:strongly_convex} with $f_i$ being replaced by $g_i$. Let us consider the inexact gradient descent (IGD) on the function $g$:
\begin{equation}\label{eq:IGM}
p^{k+1} = p^k - \theta \frac{1}{n} \sum_{i=1}^n \dgi(s_i^k),
\end{equation}
where $\theta$ is the step-size. Note that since this subsection has nothing to do with time-varying setup, to avoid heavy notation, we use the upper right corner $p^k$ instead of $p(k)$ to denote the value of $p$ at iteration $k$. In particular, we use $(p)^a$ instead of $p^a$ to denote the $a$-th power of $p$ when it may cause confusion. Let $p^*$ be the global minimum of $g$, and define
\[
r^k\triangleq \|p^k - p^*\|_\Fro\text{ for any }k=0,1,\ldots.
\]

The main lemma of this subsection is stated next; it is basically obtained by following the ideas in \cite{Devolder2013}.

\smallskip
\begin{lemma}[\textbf{The error bound on the IGD}]\label{lemma:graderror}
Suppose that
\begin{equation}\label{eq:IGM_lambda}
\sqrt{1-\frac{\theta\barmu\beta}{\beta+1}}\leq\lambda<1\quad\text{and}\quad \theta\leq\frac{1}{(1+\eta)\barL},
\end{equation}
where $\beta>0$ and $\eta>0$.  Then under Assumptions \ref{ass:smooth} and \ref{ass:strongly_convex} with $f_i$'s replaced by $g_i$'s, the tuple sequence $\{r^k, p^k; s_1^k,s_2^k,\ldots,s_n^k\}$ generated by the inexact gradient method \eqref{eq:IGM} obeys
\[\label{eq:IGM_ineq}
\begin{array}{l}
|r|^{\lambda,K}\leq2r^0
+(\lambda\sqrt{n})^{-1}\left(\sqrt{\frac{L(1+\eta)}{\barmu\eta}+\frac{\hatmu}{\barmu}\beta}\right)\sum\limits_{i=1}^n\|p-s_i\|_{\Fro}^{\lambda,K}\text{ for any }K=0,1,\ldots.
\end{array}
\]
\end{lemma}
\smallskip
\begin{proof}
By assumptions, for each $i\in\{1,2,\ldots,n\}$ and $k=0,1,\ldots$, we have
\begin{equation}\label{eq:IGM_proof1}
g_i(p^*) \geq g_i(s^k_i) + \langle\dgi(s^k_i),p^* - s^k_i\rangle + \frac{\mu_i}{2} \| p^* - s^k_i \|_\Fro^2.
\end{equation}
Through using the basic inequality $\|s^k_i - p^*\|_\Fro^2 \geq \frac{\beta}{\beta+1}\|p^k - p^*\|_\Fro^2 - \beta\|p^k - s^k_i\|_\Fro^2$ where $\beta>0$ is a tunable parameter, it follows from \eqref{eq:IGM_proof1} that
\[
\begin{array}{rcl}
g_i(p^*)
&\geq&g_i(s^k_i) +\langle\dgi(s^k_i),p^k - s^k_i\rangle +  \langle\dgi(s^k_i),p^* - p^k\rangle\\
&    &+ \frac{\mu_i}{2} \left( \frac{\beta}{\beta+1}\|p^k - p^*\|_\Fro^2 - \beta\|p^k - s^k_i\|_\Fro^2 \right)
\end{array}
\]
and therefore
\begin{equation}\label{eq:IGM_proof3}
\begin{array}{rcl}
\langle\dgi(s^k_i), p^* - p^k\rangle
&\leq& g_i(p^*) - g_i(s^k_i) - \langle\dgi(s^k_i), p^k - s^k_i\rangle\\
&    & - \frac{\mu_i\beta}{2(\beta+1)} \|p^k - p^*\|_\Fro^2 + \frac{\mu_i\beta}{2} \|s^k_i - p^k\|_\Fro^2.
\end{array}
\end{equation}
Averaging \eqref{eq:IGM_proof3} over $i$ through $1$ to $n$ gives
\begin{equation} \label{eq:IGM_lb}
\begin{array}{rcl}
&    &\frac{1}{n} \sum_{i=1}^n \langle\dgi(s^k_i), p^* - p^k\rangle\\
&\leq&g(p^*) - \frac{1}{n} \sum_{i=1}^n \left( g_i(s^k_i) + \langle\dgi(s^k_i), p^k - s^k_i\rangle - \frac{\mu_i\beta}{2} \|s^k_i - p^k\|_\Fro^2 \right)\\
&    &- \frac{\barmu\beta}{2(\beta+1)} \|p^k - p^*\|_\Fro^2.
\end{array}
\end{equation}
On the other hand, we also have that for any vector $\Delta$,
\[
\begin{array}{rcl}
g_i(p^k + \Delta)
& =  & g_i \left(s^k_i + \Delta + p^k - s^k_i \right) \\
&\leq& g_i(s^k_i) + \langle\dgi(s^k_i), \Delta  + p^k - s^k_i\rangle + \frac{L_i}{2} \|\Delta + p^k - s^k_i\|_\Fro^2 \\
&\leq& g_i(s^k_i) + \langle\dgi(s^k_i), \Delta\rangle + \langle\dgi(s^k_i),p^k - s^k_i\rangle\\
&    & + \frac{L_i(1+\eta)}{2} \|\Delta\|_\Fro^2  + \frac{L_i(1+\eta)}{2\eta} \|p^k - s^k_i\|_\Fro^2
\end{array}
\]
where $\eta>0$ is some tunable parameter, and therefore
\begin{equation}\label{eq:IGM_proof4}
\begin{array}{rcl}
- \langle\dgi(s^k_i), \Delta\rangle
&\leq& - g_i(p^k + \Delta) + g_i(s^k_i) + \langle\dgi(s^k_i), p^k - s^k_i\rangle\\
&    &+ \frac{L_i(1+\eta)}{2\eta} \|p^k - s^k_i\|_\Fro^2 + \frac{L_i(1+\eta)}{2}\|\Delta\|_\Fro^2.
\end{array}
\end{equation}
Averaging \eqref{eq:IGM_proof4} over $i$ through $1$ to $n$ gives
\begin{equation}\label{eq:IGM_ub}
\hspace{-1.5em}
\begin{array}{rcl}
- \langle\frac{1}{n} \sum_{i=1}^n \dgi(s^k_i), \Delta\rangle
& \leq & - g(p^k + \Delta)+\frac{\barL(1+\eta)}{2}\|\Delta\|_\Fro^2\\
&      &+ \frac{1}{n} \sum_{i=1}^n \left( g_i(s^k_i) + \langle\dfi(s^k_i), p^k - s^k_i\rangle + \frac{L_i(1+\eta)}{2\eta} \|p^k - s^k_i\|_\Fro^2 \right).
\end{array}
\end{equation}

Having \eqref{eq:IGM_lb} and \eqref{eq:IGM_ub} at hand, we are ready to show how $r^{k+1}$ is related to $r^k$. First, plugging $a = p^{k+1} - p^*$ and $b = p^k - p^{k+1}$ into the basic equality $\|a\|_\Fro^2=\|a+b\|_\Fro^2 - 2 \langle a, b\rangle - \|b\|_\Fro^2$ yields
\begin{equation}\label{eq:IGM_proof5}
\begin{array}{rcl}
(r^{k+1})^2
& = & (r^k)^2 - 2 \langle p^{k+1} - p^*, p^k - p^{k+1} \rangle - \|p^{k+1} - p^k\|_\Fro^2\\
&   & \text{(substituting \eqref{eq:IGM})}\\
& = & (r^k)^2 - 2 \langle p^{k+1} - p^*, \theta \frac{1}{n} \sum\limits_{i=1}^n \dgi(s_i^k) \rangle - \|p^{k+1} - p^k\|_\Fro^2\\
& = & (r^k)^2 + 2 \theta \langle \frac{1}{n} \sum\limits_{i=1}^n \dgi(s_i^k), p^* - p^{k}\rangle\\
&   & - 2 \theta \langle \frac{1}{n} \sum\limits_{i=1}^n \dgi(s_i^k), p^{k+1} - p^k \rangle - \|p^{k+1} - p^k\|_\Fro^2.\\
\end{array}
\end{equation}
Next, in \eqref{eq:IGM_proof5}, we substitute \eqref{eq:IGM_lb} for the second term, and we substitute \eqref{eq:IGM_ub} with $\Delta=p^{k+1} - p^k$ for the third term. Thus, we obtain that
\begin{equation}\label{eq:IGM_proof6}
\begin{array}{rcl}
(r^{k+1})^2
&\leq& (r^k)^2 + 2\theta \left(g(p^*) - \frac{1}{n} \sum_{i=1}^n \left( g_i(s^k_i) + \langle\dgi(s^k_i), p^k - s^k_i\rangle \right.\right.\\
&    & \left.\left. - \frac{\mu_i\beta}{2} \|s^k_i - p^k\|_\Fro^2 \right) - \frac{\barmu\beta}{2(\beta+1)} \|p^k - p^*\|_\Fro^2\right) \\
&    & + 2\theta \left( - g(p^{k+1}) + \frac{1}{n} \sum_{i=1}^n \left(  g_i(s^k_i) + \langle\dgi(s^k_i), p^k - s^k_i\rangle\right.\right.\\
&    & \left.\left.+ \frac{L_i(1+\eta)}{2\eta} \|p^k - s^k_i\|_\Fro^2 \right)+\frac{\barL(1+\eta)}{2}\|p^{k+1}-p^k\|_\Fro^2 \right) - \|p^{k+1} - p^k\|_\Fro^2\\
& =  & (r^k)^2 + 2\theta ( g(p^*) - g(p^{k+1}) ) - \frac{\theta\barmu\beta}{\beta+1} \|p^k- p^*\|_\Fro^2\\
&    & + \frac{1}{n}\sum_{i=1}^n \left(\frac{\theta L_i(1+\eta)}{\eta}+\theta\mu_i\beta\right)\|p^k - s_i^k\|_\Fro^2 - (1-\theta \barL(1+\eta))\|p^{k+1}-p^k\|_\Fro^2\\
&\leq& \left( 1 - \frac{\theta\barmu\beta}{\beta+1} \right) (r^k)^2 - 2\theta ( g(p^{k+1}) - g(p^*) )\\
&    & + \left(\frac{\theta L(1+\eta)}{\eta}+\theta\hatmu\beta\right) \frac{1}{n}\sum_{i=1}^n \|p^k - s_i^k\|_\Fro^2 - (1-\theta \barL(1+\eta))\|p^{k+1}-p^k\|_\Fro^2.
\end{array}
\end{equation}
Define $\epsilon^{k}=\frac{1}{n} \sum_{i=1}^n \|p^k - s^k_i\|_\Fro^2$. By choosing $\theta\leq\frac{1}{(1+\eta)\barL}$ such that $1-\theta \barL(1+\eta)$ in \eqref{eq:IGM_proof6} is nonnegative, we have
\begin{equation}\label{eq:IGM_proof7}
\begin{array}{l}
(r^{k+1})^2 \leq \left( 1 - \frac{\theta\barmu\beta}{\beta+1} \right) (r^k)^2 - 2\theta ( g(p^{k+1}) - g(p^*) ) + \left(\frac{\theta L(1+\eta)}{\eta}+\theta\hatmu\beta\right)\epsilon^{k}.
\end{array}
\end{equation}
Let us look into the last two terms of \eqref{eq:IGM_proof7}. Noticing that $\barmu=(1/n)\sum_{i=1}^n\mu_i$ is a strong convexity constant of $g(p)$, there are two possibilities that could happen at time $k$. Possibility A is that
\[
(r^{k+1})^2\geq \left(\frac{L(1+\eta)}{\barmu\eta}+\frac{\hatmu}{\barmu}\beta\right)\epsilon^k ,
\]
while possibility B is the opposite, namely that
\[
(r^{k+1})^2 <  \left(\frac{L(1+\eta)}{\barmu\eta}+\frac{\hatmu}{\barmu}\beta\right)\epsilon^k.
\]
If possibility A occurs, we have
\[
\begin{array}{l}
2\theta(g(p^{k+1})  - g(p^*)) \geq  \theta \barmu \|p^{k+1}-p^*\|_\Fro^2=\theta\barmu (r^{k+1})^2 \geq \left(\frac{\theta L(1+\eta)}{\eta}+\theta\hatmu\beta\right)\epsilon^k
\end{array}
\]
which together with \eqref{eq:IGM_proof7} implies
\[
(r^{k+1})^2\leq\left( 1 - \frac{\theta\barmu\beta}{\beta+1} \right) (r^k)^2.
\]
Considering both possibilities A and B, it follows that
\begin{equation}\label{eq:IGM_proof8}
\begin{array}{l}
(r^{k+1})^2\leq\max\left\{\left(1-\frac{\theta\barmu\beta}{\beta+1} \right)(r^k)^2,\left(\frac{L(1+\eta)}{\barmu\eta}+\frac{\hatmu}{\barmu}\beta\right)\epsilon^k\right\}.
\end{array}
\end{equation}
Recursively using the inequality \eqref{eq:IGM_proof8} we can see that
\begin{equation}\label{eq:IGM_proof9}
\begin{array}{ll}
&(r^{k+1})^2\leq\max
\left\{\left(1-\frac{\theta\barmu\beta}{\beta+1} \right)^{k+1}(r^0)^2,
\left(\frac{L(1+\eta)}{\barmu\eta}+\frac{\hatmu}{\barmu}\beta\right)\max\limits_{t=0,\ldots,k}\left\{\left(1-\frac{\theta\barmu\beta}{\beta+1} \right)^t\epsilon^{k-t}\right\}\right\}.
\end{array}
\end{equation}
Taking square root on both sides of \eqref{eq:IGM_proof9} gives us
\begin{equation}\label{eq:IGM_proof10}
\begin{array}{l}
r^{k+1}\leq
\left(\sqrt{1-\frac{\theta\barmu\beta}{\beta+1}}\right)^{k+1}r^0+\left(\sqrt{\frac{L(1+\eta)}{\barmu\eta}+\frac{\hatmu}{\barmu}\beta}\right)\max\limits_{t=0,\ldots,k}\left\{\left(\sqrt{1-\frac{\theta\barmu\beta}{\beta+1}} \right)^t\sqrt{\epsilon^{k-t}}\right\}.
\end{array}
\end{equation}

Choose $\lambda$ that satisfies \eqref{eq:IGM_lambda} so to have $c\triangleq(\lambda)^{-2}\left(1-\frac{\theta\barmu\beta}{\beta+1}\right)\leq1$, then from \eqref{eq:IGM_proof10} we get
\begin{equation}\label{eq:IGM_proof11}
\begin{array}{rcl}
(\lambda)^{-(k+1)}r^{k+1}
&\leq&\left(\sqrt{c}\right)^{k+1}r^0+(\lambda)^{-1}\left(\sqrt{\frac{L(1+\eta)}{\barmu\eta}+\frac{\hatmu}{\barmu}\beta}\right)\max\limits_{t=0,\ldots,k}\left\{(\lambda)^{-(k-t)}\left(\sqrt{c} \right)^t\sqrt{\epsilon^{k-t}}\right\}\\
&\leq&r^0+(\lambda)^{-1}\left(\sqrt{\frac{L(1+\eta)}{\barmu\eta}+\frac{\hatmu}{\barmu}\beta}\right)\max\limits_{t=0,\ldots,k}\left\{(\lambda)^{-t}\sqrt{\epsilon^{t}}\right\}.
\end{array}
\end{equation}
Further observing that
\[
\sqrt{\epsilon^k} = \sqrt{\frac{1}{n} \sum_{i=1}^n \|p^k - s^k_i\|_\Fro^2} \leq \frac{1}{\sqrt{n}} \sum\limits_{i=1}^n \|p^k - s^k_i\|_\Fro,
\]
and combining it with \eqref{eq:IGM_proof11}, it follows that
\begin{equation}\label{eq:IGM_proof12}
\begin{array}{l}
(\lambda)^{-(k+1)}r^{k+1}\leq r^0+(\lambda\sqrt{n})^{-1}\left(\sqrt{\frac{L(1+\eta)}{\barmu\eta}+\frac{\hatmu}{\barmu}\beta}\right)\sum\limits_{i=1}^n\max\limits_{t=0,\ldots,k}\left\{(\lambda)^{-t}\|p^t - s^t_i\|_\Fro\right\}.
\end{array}
\end{equation}
Taking $\max_{k=0,1,\ldots,K-1}\{\cdot\}$ on both sides of \eqref{eq:IGM_proof12} gives
\[\label{eq:IGM_proof13}
\begin{array}{rcl}
|r|^{\lambda,K}
&\leq&2r^0+(\lambda\sqrt{n})^{-1}\left(\sqrt{\frac{L(1+\eta)}{\barmu\eta}+\frac{\hatmu}{\barmu}\beta}\right)\sum\limits_{i=1}^n\|p-s_i\|_{\Fro}^{\lambda,K-1}\\
&\leq&2r^0+(\lambda\sqrt{n})^{-1}\left(\sqrt{\frac{L(1+\eta)}{\barmu\eta}+\frac{\hatmu}{\barmu}\beta}\right)\sum\limits_{i=1}^n\|p-s_i\|_{\Fro}^{\lambda,K}.\\
\end{array}
\]
\end{proof}
\smallskip

\subsection{The Last Arrow}

Now we prove the last arrow of our proof sketch [cf. \eqref{eq:cycle_alg1}] in the following lemma. Its establishment will use the error bound on the IGD of Lemma \ref{lemma:graderror}, as a key ingredient.

\smallskip
\begin{lemma}[\textbf{Algorithm 1: The last arrow $\cbx\rightarrow\bq$}]\label{lemma:last_arrow}
Let Assumptions \ref{ass:matrix_W}, \ref{ass:smooth}, and \ref{ass:strongly_convex} hold. In addition, assume that the stepsize $\alpha>0$ and the parameter $\lambda$ are such that
\[\label{eq:last_arrow_lambda}
\sqrt{1-\frac{\alpha\barmu\beta}{\beta+1}}\leq\lambda<1\quad\text{and}\quad\alpha\leq\frac{1}{(1+\eta)\barL},
\]
where $\beta>0$ and $\eta>0$ are some tunable parameters. Then, we have
\begin{equation}\label{eq:last_arrow_proof0}
\begin{array}{l}
\|\bq\|_\Fro^{\lambda, K} \leq \left( 1 + \frac{\sqrt{n}}{\lambda}\sqrt{\frac{L(1+\eta)}{\barmu\eta}+\frac{\hatmu}{\barmu}\beta} \right) \|\cbx\|_\Fro^{\lambda, K}+ 2\sqrt{n} \|\barx(0) - x^*\|_\Fro\text{ for any }K=0,1,\ldots.
\end{array}
\end{equation}
\end{lemma}
\smallskip
\begin{proof}
First, let us consider the evolution of $\overline{x}(k)$. Noticing that
\[
\begin{array}{rcl}
\bary(k+1)-\frac{1}{n}(\df(\bx(k+1)))^\T\one
&=&\bary(k)-\frac{1}{n}(\df(\bx(k)))^\T\one\\
& \vdots &\\
&=&\bary(0)-\frac{1}{n}(\df(\bx(0)))^\T\one\\
&=&0
\end{array}
\]
holds for all $k$, we then have that
\begin{equation}\label{eq:last_arrow_proof1}
\begin{array}{rcl}
\overline{x}(k+1)
 &=& \overline{x}(k) - \alpha\overline{y}(k)\\
 &=& \overline{x}(k) - \alpha\frac{1}{n}(\df(\bx(k)))^\T\one\\
 &=& \overline{x}(k) - \alpha\frac{1}{n}\sum\limits_{i=1}^n \nabla f_i(x_i(k)).
\end{array}
\end{equation}
Applying Lemma \ref{lemma:graderror} to the recursion relation of $\barx$, namely \eqref{eq:last_arrow_proof1}, we obtain
\begin{equation}\label{eq:last_arrow_proof2}
\hspace{-2em}\begin{array}{rcl}
\|\overline{x} - x^*\|_\Fro^{\lambda, K}
&\leq&2\|\overline{x}(0) - x^*\|_\Fro
+(\lambda\sqrt{n})^{-1}\left(\sqrt{\frac{L(1+\eta)}{\barmu\eta}+\frac{\hatmu}{\barmu}\beta}\right)\sum\limits_{i=1}^n \|\overline{x} - x_i\|_\Fro^{\lambda, K} \\
&  \leq &2\|\overline{x}(0) - x^*\|_\Fro
+(\lambda\sqrt{n})^{-1}\left(\sqrt{\frac{L(1+\eta)}{\barmu\eta}+\frac{\hatmu}{\barmu}\beta}\right)\sqrt{n}\sqrt{\sum\limits_{i=1}^n \left(\|\overline{x} - x_i\|_\Fro^{\lambda, K}\right)^2} \\
&  \leq &2\|\overline{x}(0) - x^*\|_\Fro+(\lambda)^{-1}\left(\sqrt{\frac{L(1+\eta)}{\barmu\eta}+\frac{\hatmu}{\barmu}\beta}\right)\|\cbx\|_\Fro^{\lambda, K}.
\end{array}
\end{equation}
Since
\[\label{eq:last_arrow_proof3}
\begin{array}{rcl}
\bq(k)
& = &\bx(k) - \one(\overline{x}(k))^\T + \one(\overline{x}(k))^\T  - \bx^*\\
& = &\cbx(k) + \one(\overline{x}(k)-x^*)^\T,
\end{array}
\]
it follows that
\begin{equation}\label{eq:last_arrow_proof4}
\|\bq\|_\Fro^{\lambda, K} \leq \|\cbx\|_\Fro^{\lambda, K} + \sqrt{n}\|\overline{x} - x^*\|_\Fro^{\lambda, K}.
\end{equation}
Substituting \eqref{eq:last_arrow_proof2} into \eqref{eq:last_arrow_proof4} yields \eqref{eq:last_arrow_proof0}.
\end{proof}
\smallskip

\subsection{Linear Convergence of DIGing}

We now state our main results on the convergence rates of DIGing. Our first theorem gives an explicit convergence rate for DIGing in terms of the network parameters ($B$, $n$, and $\sig$), objective parameters ($\barmu$ and $\barkappa=\frac{L}{\barmu}$), and the algorithmic step-size ($\alpha$).

\smallskip
\begin{theorem}[\textbf{Algorithm 1: Explicit geometric rate over time-varying graphs}]\label{theorem:final_bound}
Suppose that Assumptions \ref{ass:matrix_W}, \ref{ass:smooth}, and \ref{ass:strongly_convex} hold. Let
\[
\sig=\sup_{k\geq B-1}\left\{\sigmax{\Wdo_B(k)-\frac{1}{n}\one\one^\T}\right\}\quad\text{and}\quad\Cn=3\barkappa B^2\left(1+4\sqrt{n}\sqrt{\barkappa}\right).
\]
Then, for any step-size $\alpha\in\left(0,\frac{1.5(1-\sig)^2}{\barmu\Cn}\right]$, the sequence $\{\bx(k)\}$ generated by DIGing algorithm converges to the matrix $\bx^*=\one(x^*)^\T$, where $x^*$ is the unique optimal solution of problem~\eqref{eq:F} at a global R-linear (geometric) rate $O(\lambda^k)$, where the parameter $\lambda$ is given by
\[\label{eq:final_rate}
\lambda=\left\{
         \begin{array}{ll}
           \sqrt[2B]{1-\frac{\alpha\barmu}{1.5}}, &\text{ if } \alpha\in\left(0,\frac{1.5\left(\sqrt{\Cn^2+(1-\sig^2)\Cn}-\sig\Cn\right)^2}{\barmu\Cn(\Cn+1)^2}\right],\\
           \sqrt[B]{\sqrt{\frac{\alpha\barmu\Cn}{1.5}}+\sig}, &\text{ if } \alpha\in\left(\frac{1.5\left(\sqrt{\Cn^2+(1-\sig^2)\Cn}-\sig\Cn\right)^2}{\barmu\Cn(\Cn+1)^2},\frac{1.5(1-\sig)^2}{\barmu\Cn}\right].
         \end{array}
       \right.
\]
\end{theorem}
\smallskip
\begin{proof}
Let us collect all the relations/arrows at hand [cf. Lemmata \ref{lemma:first_arrow}, \ref{lemma:second_arrow}, \ref{lemma:third_arrow}, and \ref{lemma:last_arrow}]:
\[
\begin{array}{ll}
\text{i)}&\ \|\bz\|_\Fro^{\lambda, K} \leq \gamma_1 \|\bq\|_\Fro^{\lambda, K}+\omega_1 \text{ where } \gamma_1=L\left(1+\frac{1}{\lambda}\right) \text{ and } \omega_1=0;\\
\text{ii)}&\ \|\cby\|_\Fro^{\lambda, K} \leq \gamma_2 \|\bz\|_\Fro^{\lambda, K} + \omega_2 \text{ where } \gamma_2=\frac{\lambda(1-\lambda^B)}{(\lambda^B-\sig)(1-\lambda)} \text{ and } \omega_2=\frac{\lambda^B}{\lambda^B-\sig}\sum\limits_{t=1}^B\lambda^{1-t}\|\cby(t-1)\|_{\Fro};\\
\text{iii)}&\ \|\cbx\|_\Fro^{\lambda, K} \leq \gamma_3 \|\cby\|_\Fro^{\lambda, K} + \omega_3 \text{ where } \gamma_3=\frac{\alpha(1-\lambda^B)}{(\lambda^B-\sig)(1-\lambda)} \text{ and } \omega_3 = \frac{\lambda^B}{\lambda^B-\sig}\sum\limits_{t=1}^{B}\lambda^{1-t}\|\cbx(t-1)\|_\Fro;\\
\text{iv)}&\ \|\bq\|_\Fro^{\lambda, K} \leq \gamma_4 \|\cbx\|_\Fro^{\lambda, K} + \omega_4 \text{ where } \gamma_4=1 + \frac{\sqrt{n}}{\lambda}\sqrt{\frac{L(1+\eta)}{\barmu\eta}+\frac{\hatmu}{\barmu}\beta} \text{ and } \omega_4=2\sqrt{n} \|\barx(0) - x^*\|_\Fro.
\end{array}
\]
To apply the small gain theorem (Theorem \ref{theorem:small}) to get that $\|\bq\|_\Fro^{\lambda}$ is bounded, we need $\gamma_1\gamma_2\gamma_3\gamma_4<1$, that is,
\begin{equation}\label{eq:final_1}
\begin{array}{l}
\frac{\alpha L(\lambda+1)(1-\lambda^B)^2}{(\lambda^B-\sig)^2(1-\lambda)^2}\left(1 + \frac{\sqrt{n}}{\lambda}\sqrt{\frac{L(1+\eta)}{\barmu\eta}+\frac{\hatmu}{\barmu}\beta}\right)<1,
\end{array}
\end{equation}
\begin{equation}\label{eq:final_2}
\text{ where } \beta>0\quad\text{and}\quad\eta>0,
\end{equation}
along with other restrictions on parameters that appear in Lemmata \ref{lemma:second_arrow}, \ref{lemma:third_arrow}, and \ref{lemma:last_arrow}:
\begin{equation}\label{eq:final_3}
\sig < \lambda^B < 1;
\end{equation}
\begin{equation}\label{eq:final_4}
\sqrt{1-\frac{\alpha\barmu\beta}{\beta+1}}\leq\lambda<1;
\end{equation}
\begin{equation}\label{eq:final_5}
\text{and}\quad\alpha\leq\frac{1}{(1+\eta)\barL}.
\end{equation}

We next use relations \eqref{eq:final_1}--\eqref{eq:final_5} with a specific values for the parameters $\beta, \eta$ and the stepsize $\alpha$, which yields the desired result. Specifically, let $\beta=2L/\hatmu$ and $\eta=1$ in relation \eqref{eq:final_2}. By further using $0.5\leq\lambda<1$ and $(1-\lambda^B)/(1-\lambda)\leq B$, \eqref{eq:final_1} and \eqref{eq:final_5}
together yield
\begin{equation}\label{eq:final_proof1}
\alpha\leq\frac{(\lambda^B-\sig)^2}{2LB^2\left(1+4\sqrt{n}\sqrt{\barkappa}\right)}.
\end{equation}
On the other hand, since $1.5\geq 1+1/\beta$, relation \eqref{eq:final_4} implies that
\begin{equation}\label{eq:final_proof2}
\alpha\geq\frac{1.5(1-\lambda^2)}{\barmu}.
\end{equation}
Using \eqref{eq:final_proof1} and \eqref{eq:final_proof2}, it remains to show that there exists $\lambda\in(\sqrt[B]{\sig},1)$ [cf. \eqref{eq:final_3}] such that
\[\label{eq:final_proof3}
\left[\frac{1.5(1-\lambda^2)}{\barmu},\frac{(\lambda^B-\sig)^2}{2LB^2\left(1+4\sqrt{n}\sqrt{\barkappa}\right)}\right]\neq\emptyset,
\]
which is equivalent to
\begin{equation}\label{eq:final_proof4}
\left[\frac{1.5(1-\lambda^2)}{\barmu},\frac{1.5(\lambda^B-\sig)^2}{\barmu\Cn}\right]\neq\emptyset,
\end{equation}
where $\Cn=3\barkappa B^2\left(1+4\sqrt{n}\sqrt{\barkappa}\right)$. We consider a smaller interval by enlarging the left-bound of the interval in \eqref{eq:final_proof4}, i.e., we will prove that
\begin{equation}\label{eq:final_proof5}
\left[\frac{1.5(1-\lambda^{2B})}{\barmu},\frac{1.5(\lambda^B-\sig)^2}{\barmu\Cn}\right]\neq\emptyset.
\end{equation}
When $\lambda$ varies from $\sqrt[B]{\sig}$ to $1$, the left-bound of the interval in \eqref{eq:final_proof5} is monotonically decreasing from $\frac{1.5(1-\sig^2)}{\barmu}$ to $0$, while its right-bound is monotonically increasing from $0$ to $\frac{1.5(1-\sig)^2}{\barmu \Cn}$. In particular, the  relation \eqref{eq:final_proof5} is valid when $\lambda$ (as small as the current choice of all parameters can give) is given by
\[\label{eq:final_proof6}
\lambda=\sqrt[B]{\frac{\sqrt{\Cn^2+(1-\sig^2)\Cn}+\sig}{\Cn+1}}.
\]
Thus, for
\[
\alpha\in\left(0,\frac{1.5\left(\sqrt{\Cn^2+(1-\sig^2)\Cn}-\sig\Cn\right)^2}{\barmu\Cn(\Cn+1)^2}\right],
\]
we can set $\lambda=\sqrt[2B]{1-\frac{\alpha\barmu}{1.5}}$, while for
\[
\alpha\in\left(\frac{1.5\left(\sqrt{\Cn^2+(1-\sig^2)\Cn}-\sig\Cn\right)^2}{\barmu\Cn(\Cn+1)^2},\frac{1.5(1-\sig)^2}{\barmu\Cn}\right],
\]
we can use $\lambda=\sqrt[B]{\sqrt{\frac{\alpha\barmu\Cn}{1.5}}+\sig}$. The rest of the statements follow from Theorem \ref{theorem:small} and Lemma \ref{lemma:bound_R-linear}.
\end{proof}
\smallskip

Other possible choices of $\beta$, $\eta$, $\alpha$, and $\lambda$ exist and may give tighter bounds but here we only aim to give an explicit estimation on the rate.

To see how the geometric rate scales with the number of agents, we further have the following corollary.

\smallskip
\begin{corollary}[\textbf{Algorithm 1: Polynomial network scalability}]\label{corollary:scalability} Under Assumptions \ref{ass:connectedness}, \ref{ass:matrix_W_specific}, \ref{ass:smooth}, and \ref{ass:strongly_convex}, if the agents choose the step-size to be
\begin{equation}\label{eq:final_step_size}
\alpha(\tau)=\frac{3\tau^2}{128B^2 n^{4.5} L\sqrt{\barkappa}}-\frac{1.5}{\barmu}\left(\frac{\tau^2}{128B^2 n^{4.5} \barkappa^{1.5}}\right)^2,
\end{equation}
where $\tau$ is the smallest nonzero positive element of the nonnegative matrices $\Wdo(k)$ for all $k$ [cf.\ Assumption \ref{ass:matrix_W_specific}]. Then, the sequence $\{{\bf x}(k)\}$ generated by DIGing converges to the unique optimal solution ${\bf x}^*$ at a global R-linear (geometric) rate of $O\left((\lambda(\tau))^k\right)$ where
\[\label{eq:final_rate2}
\lambda(\tau)=\sqrt[B]{1-\frac{\tau^2}{128B^2 n^{4.5} \barkappa^{1.5}}}.
\]
\end{corollary}
\smallskip
\begin{proof}
Define $\varphi=1-\lambda^B$, then requiring that the interval in \eqref{eq:final_proof5} be nonempty is equivalent to showing that the following inequality has a solution:
\[\label{eq:coro_proof1}
(1-(1-\varphi)^2)\Cn\leq(1-\sig-\varphi)^2.
\]
Therefore, we can show that an achievable $\varphi$ is
\begin{equation}\label{eq:coro_proof2}
\begin{array}{rcl}
\varphi&=&\frac{2\Cn+2(1-\sig)-\sqrt{(2\Cn+2(1-\sig))^2-4(1-\sig)^2(\Cn+1)}}{2(\Cn+1)}\\
    &=&\frac{2(1-\sig)^2}{2\Cn+2(1-\sig)+\sqrt{(2\Cn+2(1-\sig))^2-4(1-\sig)^2(\Cn+1)}}\\
    &\geq&\frac{(1-\sig)^2}{2(\Cn+1)}.
\end{array}
\end{equation}
By Assumption \ref{ass:matrix_W_specific}, from Lemma 9 of reference \cite{Nedic2009_2}, we have that $\sig\leq1-\frac{\tau}{2n^2}$. Substituting $\Cn=3B^2\barkappa\left(1+4\sqrt{n}\sqrt{\barkappa}\right)$ into \eqref{eq:coro_proof2} gives us
\[\label{eq:coro_proof3}
\varphi\geq\frac{\tau^2}{8n^4(3B^2\barkappa\left(1+4\sqrt{n}\sqrt{\barkappa}\right)+1)}\geq\frac{\tau^2}{128 B^2 n^{4.5} \barkappa^{1.5}}.
\]
Thus the final rate is $\lambda=\sqrt[B]{1-\frac{\tau^2}{128 B^2 n^{4.5} \barkappa^{1.5}}}$, and a step-size to reach this rate is $\alpha=\frac{1.5(1-\lambda^{2B})}{\barmu}=\frac{3\tau^2}{128B^2 n^{4.5}L\sqrt{\barkappa}}-\frac{1.5}{\barmu}\left(\frac{\tau^2}{128B^2 n^{4.5} \barkappa^{1.5}}\right)^2$.
\end{proof}
\smallskip

Corollary \ref{corollary:scalability} explicitly shows how the linear convergence rate of the DIGing algorithm depends on the condition number $\barkappa$, time-varying graph connectivity constant $B$, and the network size $n$. To reach $\varepsilon$-accuracy, the iteration complexity under conditions of Corollary \ref{corollary:scalability} is $O\left(\tau^{-2}B^3n^{4.5}\barkappa^{1.5}\ln{\frac{1}{\varepsilon}}\right)$ which is polynomial in the number of agents $n$. Beyond the more general form of it, the advantage of Theorem \ref{theorem:final_bound} is that it explicitly depends on the parameter $\sig$ which measures the convergence speed of consensus. Indeed, Corollary \ref{corollary:scalability} uses the bound $\sig \leq 1 - \tau/(2n^2)$ from \cite{Nedic2009_2}, which may be very conservative since it applies to a rather general class of graphs. Moreover, any further advances in ``consensus theory'' deriving improved convergence bounds on consensus would immediately translate into improvements via Corollary \ref{corollary:scalability} [cf. \eqref{eq:coro_proof2}], where better bounds would immediately arise with smaller values of $\sig$.

\smallskip
\begin{corollary}[\textbf{Algorithm 1: Iteration complexity under lazy Metropolis mixing}]\label{corollary:Metropolis} Suppose Assumptions \ref{ass:connectedness}, \ref{ass:smooth}, and \ref{ass:strongly_convex} hold. Also assume that the graphs are time-invariant, undirected and connected (i.e., $B=1$). Let each $\Wdo(k)$ be a lazy Metropolis matrix, that is,
\[
\wdo_{ij}(k)=\left\{
         \begin{array}{ll}
           1/\left(2\max\{d_i(k),d_j(k)\}\right),     &\text{ if } ({j,i})\in\E,\\
           0,                              &\text{ if } ({j,i})\notin\E \text{ and } j\neq i,\\
           1-\sum_{l\in\Ni(k)}\wdo_{il}(k),&\text{ if } j=i.
         \end{array}
       \right.
\]
Then, with the agents choosing the step-size $\alpha(\frac{2}{71})$ [cf.\ \eqref{eq:final_step_size} with $\tau=\frac{2}{71}$ and $B=1$], the sequence $\{{\bf x}(k)\}$ generated by DIGing converges to the unique optimal solution ${\bf x}^*$ at a global R-linear (geometric) rate of $O(\lambda^k)$, where $\lambda=1-\frac{1}{161312n^{4.5} \barkappa^{1.5}}$. In particular, the number of iterations needed to reach $\varepsilon$-accuracy is $O\left(n^{4.5}\barkappa^{1.5}\ln{\frac{1}{\varepsilon}}\right)$.
\end{corollary}
\smallskip

We omit the proof for Corollary \ref{corollary:Metropolis} since it is essentially the same as the proof for Corollary \ref{corollary:scalability} in addition to which we further used $\delta =1-\frac{1}{71n^2}$ from Lemma 2.2 of reference \cite{Olshevsky2014}.

\begin{remark}[\textbf{Conservatism on the scalability of consensus speed $\delta$}]\label{remark:bound_on_delta}
\mr{In certain cases, the bound $\delta \leq 1 - \tau/(2n^2)$ can be  conservative. In the most dramatic case,  for the (fixed) complete graph this bound tells us that $1-\delta$ is bounded below by $1/n^3$, whereas it is not too hard to see that for the complete graph $1-\delta$ is actually bounded away from zero by a constant.}

\mr{In general, however, this bound cannot be improved, in the sense that there are graphs for which it is essentially tight. For example, on the (fixed) line or ring graph, the bound gives that $1-\delta$ is lower bounded by $1/n^2$, which is the correct scaling up to a constant,  as can be seen by observing that it can take two random walks at least $\Omega(n^2)$ steps to intersect on these graphs, and thus the spectral gap has to be at least that much (for more details on making such arguments rigorous, we refer the reader to the introductory chapter of  \cite{lindvall2002lectures}).}

\mr{In general, it is not possible to give a nonconservative bound on $\delta$ in terms of combinatorial features of the underlying graphs as well as the number of nodes, even for fixed matrices. Such bounds have been explored at great length within Markov chain literature, see e.g., the seminal paper \cite{diaconis1991geometric} or the monograph \cite{levin2009markov}, resulting in many different techniques, each giving accurate bounds on some graphs, and others.}

\mr{Nevertheless, for many sequences when the graph is {\em fixed}, the quantity $\delta$ can be bounded  accurately. The key tool is a connection between $\delta$ and the average hitting time (Theorem 11.11 of \cite{levin2009markov}) which in turn can be bounded in terms of graph resistance (see \cite{chandra1996electrical}). Using this connection, the following scalings may be obtained; we omit the details.}

\begin{itemize} \item \mr{On the complete graph, $\delta \leq 1 - \Omega(1)$.}
	
\item \mr{On the line and ring graphs, we have $\delta \leq 1 - \Omega(1/n^2)$.}
	
\item \mr{On the 2D grid and on the complete binary tree, $\delta \leq 1 - \Omega(1/(n \log n))$.}


\item \mr{On any regular graph, $\delta \leq 1-\Omega(1/n^2)$ as a consequence of the hitting time bounds of \cite{coppersmith1996random}.}
\item \mr{On the star graph and on the two-star graph (defined to be two stars on $n/2$ nodes with a connection between the centers), $\delta \leq 1-\Omega(1/n^2)$.}	

	
\end{itemize}
\end{remark}

\begin{remark}[\textbf{Adapt-then-Combine strategy and acceleration}]\label{remark:DIGing-ATC}
Our analytical framework also applies to Aug-DGM of \cite{Xu2015} when its step-size matrix $\bD$ is set to $\alpha I$ (see Subsection \ref{sec:Xu2015}). We also find that when the graph is well connected, the ATC strategy employed in Aug-DGM can improve the linear rate. But due to space limit, 
we omit any detailed discussion on this aspect. 
In the following design of a variant of DIGing for directed graphs, we also partially employed the ATC strategy to accelerate the convergence. \mr{Numerical experiments demonstrating the faster convergence of the ATC variant of DIGing (abbreviated as DIGing-ATC) are also provided in Section \ref{sec:num_exp}.}
\end{remark}
\begin{remark}[\textbf{Parameter selection}]
\mr{Allowing the agents to use uncoordinated (different) step-sizes (all satisfying an upper bound) is possible; see our very recent subsequent work \cite{nedic2016geometrically} as well as the related papers \cite{Xu2015,lu2016geometrical}. To keep the analysis in this paper concise, we do assume all agents use the same step-size $\alpha$. }

\mr{The step size selection, as instructed by the theorems and corollaries in this paper, requires the agents to agree on a few parameters through preprocessing. It suffices for agents to know a common upper bound on the number of nodes in the network $n$ and on the connectivity constant $B$. The other parameters used in step-size selection are  (the lower bound of) $\barmu$ and (the upper bound of) $L$. Since minima (and maxima) can be easily computed in a time-varying $B$-connected network in $O(nB)$ rounds by an algorithm wherein each node repeatedly takes minima (or maxima) of in-neighbors, the amount of pre-processing in computing $\barmu,L$ is essentially negligible compared to the worst-case running time of the protocol.}
\end{remark}

\begin{remark}{\textbf{Random link activation}} \mr{One important class of time-varying graphs comes from {\em random link activation}. Consider a fixed connected graph $\Gra=(\V,\E)$; the graph $\Gra(t)$ at time $t$ comes from having each link of $\Gra$ be present independently with some fixed (constant) probability. Fix $t$; it is easy to see that under random link activation the sequence $\Gra(1), \ldots, \Gra(t)$ (that is, the graph sequence up to time $t$) will be $B$-connected with the probability $1-1/{\rm poly}(t)$ if we choose $B=O (\log(nt))$.  Applying the results of this paper, for each $t$ can get a bound which holds with high probability and is geometric in $t/\log (nt)$.}

\mr{It is possible to improve on this and obtain a geometric rate in $t$ {\em with probability one}. We sketch how this can be done next. The constant $B$ should be chosen so that the graph obtained by taking all the edges that appear over $B$ steps is connected with probability at least $1/2$. This allows us to choose $B$ to be a constant independent of $n$ and $t$. We then apply the arguments of this paper to the quantities $E[||\bq(k)||], E[||\bz(k)||], E[||\overline{\by}(k)||], E[[\overline{\bx}(k)||]$ and use the small gain theorem to obtain that all of these quantities converge to zero at a geometric rate. Markov's inequality followed by the Borel-Cantelli lemma then yields that asymptotically $\bx(k)$ converges to ${\bf x}^*$ at a geometric rate with probability one . }
\end{remark}

\section{Distributed Optimization over Directed Graphs}\label{sec:dir_opt}

We now focus our attention to directed graphs. For such graphs we want to design a distributed algorithm that
can work with mixing matrices that need not be doubly stochastic. To do so, we employ the idea of push-sum protocol which relaxes the requirement of doubly stochastic mixing matrices to column stochastic matrices to achieve average consensus. We then introduce our algorithm that uses push-sum protocol for tracking the gradient average in time. The resulting algorithm is termed Push-DIGing (Algorithm 2), which we analyze later on in Section \ref{sec:conv_analysis2}.

\subsection{Motivation}
Suppose there are $n$ agents that can communicate over a static strongly connected directed graph.
Each agent $i$ initially holds row $i$ of $\bx(0)\in\R^{n\times p}$ and would like to calculate the average $\frac{1}{n}\one^\T\bx(0)$. To do so, one possible decentralized approach is to construct a doubly stochastic matrix $\Wdo$ and perform updates $\bx(k+1)=\Wdo\bx(k)$ starting from $\bx(0)$. However, in a general directed graph, construction of a doubly stochastic matrix needs a weight balancing procedure and it is costly \cite{Gharesifard2012}. This becomes even less realistic when the graph is time-varying, as the maintenance of a doubly stochastic matrix sequence needs a real time weight balancing.

Instead, if every agent knows its out-degree, it is possible for the agents to construct a column stochastic matrix $\Wco$ and perform the following recursions with initialization $\bu(0)=\bx(0)$ and $\bv(0)=\one$ to achieve the average (push-sum protocol \cite{Kempe2003}):
\[
\begin{array}{rl}
\bu\text{-update: }&\bu(k+1)=\Wco\bu(k);\\
\bv\text{-update: }&\bv(k+1)=\Wco\bv(k);\ \bV(k+1)=\dia{\bv(k+1)};\\
\bx\text{-update: }&\bx(k+1)=(\bV(k+1))^{-1}\bu(k+1).
\end{array}
\]
Intuitively, noticing that $\baru(k+1)=\baru(k)$ ($\bu(k)$ is sum preserving), rows of $\bu(k)$ are heading towards scaled averages with uneven scaling ratios across the vertices caused by the non-double stochasticity of $\Wco$ (the ratios are actually the elements of a right eigenvector of $\Wco$ corresponding to the eigenvalue $1$), while $\bV(k)$ is recording the ratios. By applying the recorded ratio inverse $(\bV(k))^{-1}$ on $\bu(k)$, the algorithm recovers the unscaled average of the rows in $\bx(k)$.

\subsection{The Push-DIGing Algorithm}\label{sec:algo_dev2}

Next we formally state Push-DIGing in \emph{Algorithm 2}.

\smallskip
\begin{center}
  {\textbf{Algorithm 2: Push-DIGing}}

  \smallskip
    \begin{tabular}{l}
    \hline
    \emph{  } Choose step-size $\alpha>0$ and pick any $\bx(0)=\bu(0)\in \R^{n \times p}$;\\
    \emph{  } Initialize $\by(0)=\df(\bx(0))$, $\bv(0)=\one\in\R^n$, and $\bV(0)=\dia{\bv(0)}$;\\
    \emph{  } \textbf{for} $k=0,1,\ldots$ \textbf{do}\\
    \qquad$\bu(k+1)=\Wco(k)\left(\bu(k)-\alpha\by(k)\right)$;\\
    \qquad$\bv(k+1)=\Wco(k)\bv(k)$; $\bV(k+1)=\dia{\bv(k+1)}$;\\
    \qquad$\bx(k+1)=(\bV(k+1))^{-1}\bu(k+1)$;\\
    \qquad$\by(k+1)=\Wco(k)\by(k)+\df(\bx(k+1))-\df(\bx(k))$;\\
    \emph{  } \textbf{end for}\\
    \hline
    \end{tabular}
\end{center}
\smallskip

Looking at the individual agents, the initialization of Push-DIGing uses an arbitrary $x_i(0)=u_i(0)\in\R^p$, and sets $y_i(0)=\nabla f_i(x_i(0))$ and $v_i(0)=1$ for $i=1,\ldots,n$. Then, at each iteration $k$, every agent $i$ performs updates, as follows:
\[
\begin{array}{l}
u_i(k+1)=\wco_{ii}(k)(u_i(k)-\alpha y_i(k))+\sum_{j\in\Niin(k)}\wco_{ij}(k)(u_j(k)-\alpha y_j(k)),\\
v_i(k+1)=\wco_{ii}(k)v_i(k)+\sum_{j\in\Niin(k)}\wco_{ij}(k)v_j(k),\\
x_i(k+1)=u_i(k+1)/v_i(k+1),\\
y_i(k+1)=\wco_{ii}(k)y_i(k)+\sum_{j\in\Niin(k)}\wco_{ij}(k)y_j(k)+\dfi(x_i(k+1))-\dfi(x_i(k)),\\
\end{array}
\]
where $\Niin(k)$ is the set of agents that can send information to agent $i$ (in-neighbors of agent $i$) at time $k$, while $\Niout(k)$ is the set of agents that can receive the information from agent $i$ (in-neighbors of agent $i$) at time $k$. (Formal definition of the sets $\Niin(k)$ and $\Niout(k)$ is deferred to Section~\ref{sec:conv_analysis2}). At every iteration $k$, each agent $i$ sends its $u_i(k)-\alpha y_i(k)$, $y_i(k)$, and $v_i(k)$ all scaled by $C_{ij}(k)$ to each of its out-neighbors $\Niout(k)$, and receives the corresponding messages from its in-neighbors $\Niin(k)$. Then, each agent $i$ updates its own $u_i(k+1)$ by summing its own $\wco_{ii}(k)\left(u_i(k)-\alpha y_i(k)\right)$ and the received $\wco_{ij}(k)\left(u_j(k)-\alpha y_j(k)\right)$ from its in-neighbors $\Niin(k)$; a similar strategy applies to the update of $v_i(k+1)$; then $x_i(k+1)$ is given by scaling $u_i(k+1)$ with $(v_i(k+1))^{-1}$; finally each agent $i$ updates its own $y_i(k+1)$ by summing its own $\wco_{ii}y_i(k)$ and the received $\wco_{ij}(k)y_j(k)$ from its in-neighbors $\Niin(k)$, and accumulating its current local gradient $\dfi(x_i(k+1))$ and subtracting its previous local gradient $\dfi(x_i(k))$ (in order to filter in only the new information contained in the most recent gradient). Unlike DIGing for undirected graphs in which each agent scales the received variables ($x_i(k)$ and $y_i(k)$) and then sums them up, in Push-DIGing, the variables ($u_i(k)-\alpha y_i(k)$, $y_i(k)$, and $v_i(k)$) are scaled before being sent out. This is due to the fact that, over directed graphs, usually a scaling weight $\Wco_{ij}$ can only be conveniently determined by the out-degree information of agent $j$ which is not available to agent $i$.

\section{Convergence Analysis for Push-DIGing}\label{sec:conv_analysis2}

In this section we conduct the convergence analysis for Push-DIGing over time-varying directed graphs. Consider a time-varying graph sequence $\{\Grad(0),\Grad(1),\ldots\}$. Every graph instance $\Grad(k)$ consists of a static set of agents $\V=\{1,2,\ldots,n\}$ and a set $\A(k)$ of time-varying arcs. An arc $(\overrightarrow{j,i})\in\A(k)$ indicates that agent $j$ can send information to agent $i$ at time (iteration) $k$. The set of in- and out-neighbors of agent $i$ at time $k$ are defined as $\Niin(k)=\left\{j\big|(\overrightarrow{j,i})\in\A(k)\right\}$ and $\Niout(k)=\left\{j\big|(\overrightarrow{i,j})\in\A(k)\right\}$, respectively.

We make the following two assumptions for the setup of time-varying directed graphs.

\smallskip
\mr{\begin{assumption}[\textbf{$\tBG$-strongly connected graph sequence}]\label{ass:connectedness2} The time-varying directed graph sequence $\Grad(k)$ is $\tBG$-strongly connected. Specifically, there exists an integer $\tBG > 0$ such that for any $t=0,1,\ldots$, the directed graph
\[\Grad_{\tBG}(t\tBG)\triangleq\left\{\V,\bigcup_{\ell=t\tBG}^{(t+1)\tBG-1}\A(\ell)\right\}\] 
is strongly connected.
\end{assumption}}
\smallskip

\mr{Let us denote \[\BG=2\tBG-1.\] 
Note that Assumption~\ref{ass:connectedness2} 
implies that $\Grad_{\BG}(k)$ is strongly connected for all $k=0,1,\ldots$.}

\begin{assumption}[\textbf{Mixing matrix sequence $\{\Wco(k)\}$}]\label{ass:matrix_C}
For any $k=0,1,\ldots$, the mixing matrix $\Wco(k)=[\wco_{ij}(k)]\in\R^{n\times n}$ is given by
\[
\wco_{ij}(k)=\frac{1}{d_j^{\mathrm{out}}(k)+1}\text{ if }(\overrightarrow{j,i}) \in \A(k), \text{ and otherwise } \wco_{ij}(k)=0,
\]
where $d_j^{\mathrm{out}}(k)=|\Njout(k)|$ is the out-degree of agent $j$ at time $k$.
\end{assumption}
\smallskip

Assumption \ref{ass:connectedness2} has been used in distributed optimization over time-varying directed graphs \cite{Nedic2015}. Similar to the case of undirected graphs, we may have other options for the weights 
$\wco_{ij}(k)$. In the existing literature on push-sum consensus protocol, the best understood are the matrices relying on the out-degree information (as in Assumption~\ref{ass:matrix_C}), which we use in establishing the bound on convergence rate of push-sum (see Lemma \ref{lemma:mixing_contraction2} and its proof). Generalizations of Assumption \ref{ass:matrix_C} may be of their own interest 
for the push-sum consensus algorithm, but that is beyond the scope of this paper.

A little algebra shows that the recursion relation of Push-DIGing is equivalent to
\begin{equation}\label{eq:push-eq}
\begin{array}{l}
\bv(k+1)=\Wco(k)\bv(k), \bV(k+1)=\dia{\bv(k+1)},\\
\bx(k+1)=\tbR(k)\left(\bx(k)-\alpha\bh(k)\right),\\
\bh(k+1)=\tbR(k)\bh(k)+(\bV(k+1))^{-1}\left(\df(\bx(k+1))-\df(\bx(k))\right),
\end{array}
\end{equation}
where $\tbR(k)=(\bV(k+1))^{-1}\Wco(k)(\bV(k))$ and $\bh(k)=(\bV(k))^{-1}\by(k)$. We note here that, under Assumptions~\ref{ass:connectedness2} and~\ref{ass:matrix_C}, it can be seen that each matrix $\bV(k)$ is invertible, and that
\begin{equation}\label{eq:norm_invV}
\|\bV^{-1}\|_{\max}^{1}\triangleq\sup_{k\geq0}\|(\bV(k))^{-1}\|_{\max}\leq n^{n\BG},
\end{equation}
where $\BG$ is the graph connectivity constant defined in Assumption \ref{ass:connectedness2}. The preceding relation follows from Corollary 2(b) of \cite{Nedic2015} (here we borrow the notation for the special norm defined in \eqref{eq:norm}). Also, we note that $\tbR(k)$ is actually a row stochastic matrix, i.e., every row of $\tbR(k)$ sums to $1$ (see Lemma 4 of\cite{Nedic2014}).

In what follows, we will use following notation
\[
\Wco_b(k)\triangleq\Wco(k)\Wco(k-1)\cdots\Wco(k+1-b)
\]
for any $k=0,1,\ldots$ and $b=0,1,\ldots$ with the convention that $\Wco_b(k)=I$ for any needed $k<0$ and $\Wco_0(k)=I$ for any $k$. The same notation rule applies to $\tbR(k)$. In the sequel, we will give an upper bound on a norm of $(I-\one\one^\T/n)\tbR_B(k)$, as provided in the next lemma. This lemma comes from the properties of push-sum protocol and can be obtained from references \cite{Nedic2008,Nedic2010}.

\begin{lemma}[\textbf{$B$-step consensus contraction}]\label{lemma:mixing_contraction2}
Under Assumptions \ref{ass:connectedness2} and \ref{ass:matrix_C}, let $B$ be an integer satisfying $B\geq\BG$ and such that
\[\label{eq:sig2}
\delta\triangleq Q_1(1-\ttau^{n\BG})^{\frac{B-1}{n\BG}} <1,
\]
where $Q_1$ and $\ttau$ are given by
\begin{equation}\label{eq:Q_1}
Q_1\triangleq2n\frac{1+\ttau^{-n\BG}}{1-\ttau^{n\BG}}\quad\text{and}\quad\ttau\triangleq\frac{1}{n^{2+n\BG}}.
\end{equation}
Then, for any $k=B-1,B,\ldots$ and any matrix $\bb$ with appropriate dimensions, if $\ba = \tbR_B(k) \bb$, i.e., $\ba =  (\bV(k+1))^{-1}\Wco_B(k)\bV(k+1-B) \bb$, we have $\|\ba\|_{\emph{\Lap}} \leq \sig\|\bb\|_{\emph{\Lap}}$.
\end{lemma}
\smallskip

\begin{proof}
It has been shown in Lemma 4 of \cite{Nedic2015} that $\tbR(k)=(\bV(k+1))^{-1}\Wco(k)\bV(k)$ is a row stochastic matrix with entries being $\tilde{R}_{ij}(k)=C_{ij}(k)v_j(k)/v_i(k+1)$, where $v_i(k)$ denotes the $i$-th entry of the vector $\bv(k)$. Denote $\tbd(k,B)\triangleq(\one^\T\tbR_B(k))^\T/n$ which is a stochastic vector. Then it follows that
\[
\begin{array}{rcl}\label{eq:contraction2_proof1}
\|\ba\|_\Lap&=&\|(I-\frac{1}{n}\one\one^\T)\tbR_B(k)\bb\|_\Fro\\
            &=&\|(\tbR_B(k)-\one(\tbd(k,B))^\T)(I-\frac{1}{n}\one\one^\T)\bb\|_\Fro\\
            &\leq&\|\tbR_B(k)-\one(\tbd(k,B))^\T\|_2\|\bb\|_\Lap.
\end{array}
\]

Now we can focus on establishment of an upper bound for $\|\tbR_B(k)-\one(\tbd(k,B))^\T\|_2$. The first step is to analyse the constituent of $\tilde{R}_{ij}(k)=C_{ij}(k)v_j(k)/v_i(k+1)$. From Corollary 2(b) of reference \cite{Nedic2014}, we know that $v_j(k)\geq\frac{1}{n^{n\BG}}$. Also, obviously we have $1/v_j(k)\geq1/n$, and for any $(\overrightarrow{j,i})\in\A(k)$, $C_{ij}(k)\geq1/n$. Based on those bounds on $C_{ij}(k)$, $v_j(k)$, and $1/v_i(k+1)$, we can have the following equivalent arc-utilization lower bound
\[
\tilde{R}_{ij}(k)\geq\ttau\triangleq\frac{1}{n^{2+n\BG}} \text{ if } (\overrightarrow{j,i})\in\A(k).
\]
Thus, we have
\begin{equation}
\begin{array}{rcl}\label{eq:contraction2_proof2}
\|\tbR_B(k)-\one(\tbd(k,B))^\T\|_2
&\leq&n\|\tbR_B(k)-\one(\tbd(k,B))^\T\|_{\max}\\
&\leq&2n\frac{1+\ttau^{-n\BG}}{1-\ttau^{n\BG}}\left(1-\ttau^{n\BG}\right)^{\frac{B-1}{n\BG}},
\end{array}
\end{equation}
where the second inequality comes from Lemma 5 of reference \cite{Nedic2010}. By \eqref{eq:contraction2_proof1} and \eqref{eq:contraction2_proof2}, we get $\|\ba\|_\Lap\leq Q_1\left(1-\ttau^{n\BG}\right)^{\frac{B-1}{n\BG}}\|\bb\|_\Lap$ with $Q_1\triangleq2n\frac{1+\ttau^{-n\BG}}{1-\ttau^{n\BG}}$.
\end{proof}
\smallskip

The convergence analysis of Push-DIGing will be based on analogous recursions illustrated in \eqref{eq:push-eq}. Similar to the proof in Section \ref{sec:conv_analysis}, we will follow the proof sketch of the small gain theorem around the cycle:
\begin{equation}\label{eq:cycle_alg3}
\text{Algorithm 2: }\bq \rightarrow \bz \rightarrow \cbh \rightarrow \cbx \rightarrow \bq.
\end{equation}

\begin{remark}[Consensus in the parameters]
In consensus-based algorithms for optimization over directed graphs, it is difficult to construct monotonically decreasing Lyapunov functions for convergence analysis \cite{Olshevsky2007} due to the presence of asymmetric operators in the iterations (arising from the asymmetric weight matrices). Also, to deal with time-varying graphs, one has to resort to time-varying or ergodic metrics \cite{Touri2012} which are not easy to construct when the consensus protocol is combined with an optimization algorithm. In this situation, conventional approaches that heavily rely on every step contraction for proving Q-linear convergence are usually inapplicable. However, by defining a special metric and utilizing the small gain theorem, we manage to conveniently analyze the introduced algorithms and establish their linear convergence rates, but without relying on the monotonic decay of any Lyapunov function associated with the recursion.
\end{remark}

\subsection{The Establishment of Each Arrow}

Noticing that that Lemma \ref{lemma:first_arrow} is a simple consequence of Assumption \ref{ass:smooth}, so it also holds for Push-DIGing. For the sake of reference convenience, we restate it as follows without proof.
\smallskip
\begin{lemma}[\textbf{Algorithm 2: The first arrow $\bq\rightarrow \bz$}]\label{lemma:first_arrow2}
Under Assumption \ref{ass:smooth}, we have that for all $K=0,1,\ldots$ and any $\lambda\in (0,1)$,
\[
\|\bz\|_\Fro^{\lambda, K} \leq L \left( 1 + \frac{1}{\lambda} \right) \|\bq\|_\Fro^{\lambda, K}.
\]
\end{lemma}
\smallskip

The next two lemmata are provided by doing almost identical arguments as those for Lemmata \ref{lemma:second_arrow} and \ref{lemma:third_arrow}: indeed, by noticing the similarity between the equivalent recursion of Push-DIGing \eqref{eq:push-eq} and the recursion of DIGing (see Algorithm 1), similar bounds to those in Lemmata \ref{lemma:second_arrow} and \ref{lemma:third_arrow} should be obtained by an application of  Lemma~\ref{lemma:mixing_contraction2}, which shows how multiplication by a row stochastic matrix $\tbR(k)$ shrinks the distance to the consensus subspace.

\smallskip
\begin{lemma}[\textbf{Algorithm 2: The second arrow $\bz\rightarrow\cbh$}]\label{lemma:second_arrow2}
Let Assumptions \ref{ass:connectedness2} and \ref{ass:matrix_C} hold, and let $\lambda$ be such that $\sig < \lambda^B < 1$, where $B$ is the constant provided in Lemma \ref{lemma:mixing_contraction2}. Then, we have
\begin{equation}
\begin{array}{c}\label{eq:second_arrow_alg3}
\|\cbh\|_{\Fro}^{\lambda,K}\leq\frac{Q_1\|\bV^{-1}\|_{\max}^1\lambda(1-\lambda^B)}{(\lambda^B-\sig)(1-\lambda)}\|\bz\|_\Fro^{\lambda,K}
+\frac{\lambda^B}{\lambda^B-\sig}\sum\limits_{t=1}^B\lambda^{1-t}\|\cbh(t-1)\|_{\Fro}\ \text{ for all }\ K=0,1,\ldots,
\end{array}
\end{equation}
where $Q_1$ and $\|\bV^{-1}\|_{\max}^{1}$ are the constants defined by \eqref{eq:Q_1} and \eqref{eq:norm_invV}, respectively.
\end{lemma}
\smallskip
\begin{proof}
The equivalent recursion of Push-DIGing involving $\bh$ and $\bz$ is
\begin{equation}\label{eq:second_arrow_proof0_alg3}
\bh(k+1)=\tbR(k)\bh(k)+(\bV(k+1))^{-1}\bz(k+1).
\end{equation}
From \eqref{eq:second_arrow_proof0_alg3}, using Lemma \ref{lemma:mixing_contraction2}, for all $k \geq B-1$, it follows that
\begin{equation}
\hspace{-1em}\begin{array}{rcl}\label{eq:second_arrow_proof1_prime_alg3}
\|\cbh(k+1)\|_\Fro
&  =   &\|\bh(k+1)\|_\Lap\\
& \leq &\left\|\tbR_B(k)\bh(k+1-B)\right\|_{\Lap} + \left\|\tbR_{B-1}(k)(\bV(k+2-B))^{-1}\bz(k+2-B)\right\|_{\Lap}\\
&      &+ \cdots + \|\tbR_1(k)(\bV(k))^{-1}\bz(k)\|_{\Lap} + \|\tbR_0(k)(\bV(k+1))^{-1}\bz(k+1)\|_{\Lap}\\
& \leq &\sig\|\cbh(k+1-B)\|_\Fro + Q_1\sum\limits_{t=1}^B\|(\bV(k+2-t))^{-1}\bz(k+2-t)\|_\Fro\\
& \leq &\sig\|\cbh(k+1-B)\|_\Fro + Q_1\|\bV^{-1}\|_{\max}^1\sum\limits_{t=1}^B\|\bz(k+2-t)\|_\Fro,
\end{array}
\end{equation}
where $Q_1$ and $\|\bV^{-1}\|_{\max}^{1}$ are the constants defined in \eqref{eq:Q_1} and \eqref{eq:norm_invV}, respectively.

Noticing the similarity between \eqref{eq:second_arrow_proof1_prime_alg3} and \eqref{eq:second_arrow_proof1_prime}, by the same argument as we have applied in the proof of Lemma \ref{lemma:second_arrow} starting from \eqref{eq:second_arrow_proof1}, we can obtain \eqref{eq:second_arrow_alg3}.
\end{proof}
\smallskip

\begin{lemma}[\textbf{Algorithm 2: The third arrow $\cbh\rightarrow\cbx$}]\label{lemma:third_arrow2}
Let Assumptions \ref{ass:connectedness2} and \ref{ass:matrix_C} hold, and let $\lambda$ be such that $\sig < \lambda^B < 1$, where $B$ is the constant provided in Lemma \ref{lemma:mixing_contraction2}. Then, we have
\begin{equation}\label{eq:third_arrow_alg3}
\hspace{-1em}\|\cbx\|_\Fro^{\lambda, K}
\leq\frac{\alpha}{\lambda^B-\sig}\left(\sig+Q_1\frac{1-\lambda^{B-1}}{1-\lambda}\right) \|\cbh\|_\Fro^{\lambda, K}  + \frac{\lambda^B}{\lambda^B-\sig}\sum\limits_{t=1}^{B}\lambda^{1-t}\|\cbx(t-1)\|_\Fro
\end{equation}
for all $K=0,1,\ldots$, where $Q_1$ is the constant as introduced in Lemma \ref{lemma:second_arrow2} (see \eqref{eq:Q_1}).
\end{lemma}
\smallskip

The equivalent recursions of Push-DIGing involving $\bx$ and $\bh$ is
\begin{equation}\label{eq:third_arrow_proof0_alg3}
\bx(k+1)=\tbR(k)\left(\bx(k)-\alpha\bh(k)\right).
\end{equation}
Noticing the similarity between \eqref{eq:third_arrow_proof0_alg3} and \eqref{eq:second_arrow_proof0_alg3}, with almost identical argument as that illustrated in the proof of Lemma \ref{lemma:second_arrow2}, we can get \eqref{eq:third_arrow_alg3}.

Similar to the proof of Lemma \ref{lemma:last_arrow}, Lemma \ref{lemma:graderror} also serves as a key ingredient in establishing the last arrow of the proof sketch for Push-DIGing. We state it in the form of a lemma as follows.

\smallskip
\begin{lemma}[\textbf{Algorithm 2: The last arrow $\cbx\rightarrow\bq$}]\label{lemma:last_arrow_alg3}
Let Assumptions \ref{ass:smooth},  \ref{ass:strongly_convex}, and \ref{ass:matrix_C} hold. Also, assume that
\[\label{eq:last_arrow_lambda_alg3}
\sqrt{1-\frac{\alpha\barmu\beta}{\beta+1}}\leq\lambda<1\quad\text{and}\quad\alpha\leq\frac{1}{(1+\eta)\barL}
\]
where $\beta>0$ and $\eta>0$ are some tunable parameters. Then, we have that for all $K=0,1,\ldots$,
\begin{equation}\label{eq:last_arrow_proof0_alg3}
\begin{array}{l}
\|\bq\|_\Fro^{\lambda, K} \leq (1+\sqrt{n})\left( 1 + \frac{\sqrt{n}}{\lambda}\sqrt{\frac{L(1+\eta)}{\barmu\eta}+\frac{\hatmu}{\barmu}\beta} \right) \|\cbx\|_\Fro^{\lambda, K}+2\sqrt{n} \|\barx(0) - x^*\|_\Fro.
\end{array}
\end{equation}
\end{lemma}
\smallskip
\begin{proof}
First, by the same argument as in the proof of Lemma \ref{lemma:last_arrow}, we have
\begin{equation}\label{eq:last_arrow_proof1_alg3}
\begin{array}{c}
\baru(k+1)=\baru(k) - \alpha\frac{1}{n}\sum\limits_{i=1}^n f_i(x_i(k)).
\end{array}
\end{equation}

Applying Lemma \ref{lemma:graderror} to the recursion relation of $\baru$, namely \eqref{eq:last_arrow_proof1_alg3}, we obtain
\begin{equation}\label{eq:last_arrow_proof2_alg3}
\begin{array}{c}
\|\baru - x^*\|_\Fro^{\lambda, K}
\leq2\|\baru(0) - x^*\|_\Fro
+(\lambda\sqrt{n})^{-1}\left(\sqrt{\frac{L(1+\eta)}{\barmu\eta}+\frac{\hatmu}{\barmu}\beta}\right)\sum\limits_{i=1}^n \|\baru - x_i\|_\Fro^{\lambda, K}.
\end{array}
\end{equation}
Let us look into the summation in the last term of \eqref{eq:last_arrow_proof2_alg3}. Since $\bu(k)=\bV(k)\bx(k)$, it follows that
\begin{equation}\label{eq:last_arrow_proof2_2_alg3}
\begin{array}{rl}
\sum\limits_{i=1}^n \|\baru - x_i\|_\Fro^{\lambda, K}
  =  &\sum\limits_{i=1}^n \|(\baru - \barx) + (\barx - x_i)\|_\Fro^{\lambda, K}\\
\leq &n \|\baru - \barx\|_\Fro^{\lambda, K} + \sum\limits_{i=1}^n\|\barx - x_i\|_\Fro^{\lambda, K}\\
\leq &n\|\frac{1}{n}\bx^\T\bV\one - \frac{1}{n}\bx^\T\one\|_\Fro^{\lambda, K}+\sqrt{n}\|\cbx\|_\Fro^{\lambda, K}\\
  =  &\|(\one-\bv)^\T\bx\|_\Fro^{\lambda, K}+\sqrt{n}\|\cbx\|_\Fro^{\lambda, K}.
\end{array}
\end{equation}
Thus, by \eqref{eq:last_arrow_proof2_alg3} and \eqref{eq:last_arrow_proof2_2_alg3} we have
\[\label{eq:last_arrow_proof2_3_alg3}
\begin{array}{rcl}
\|\baru - x^*\|_\Fro^{\lambda, K}
&\leq&2\|\barx(0) - x^*\|_\Fro+(\lambda)^{-1}\left(\sqrt{\frac{L(1+\eta)}{\barmu\eta}+\frac{\hatmu}{\barmu}\beta}\right)\|\cbx\|_\Fro^{\lambda, K}\\
&    &+(\lambda\sqrt{n})^{-1}\left(\sqrt{\frac{L(1+\eta)}{\barmu\eta}+\frac{\hatmu}{\barmu}\beta}\right)\|(\one-\bv)^\T\bx\|_\Fro^{\lambda, K}
\end{array}
\]
Since
\[\label{eq:last_arrow_proof3_alg3}
\begin{array}{rcl}
\bq(k)
& = & \bx(k) - \one(\barx(k))^\T + \one(\barx(k))^\T - \one(\baru(k))^\T + \one(\baru(k))^\T  - \bx^*\\
& = & \cbx(k) + \frac{1}{n}\one(\one-\bv(k))^\T\bx(k)+\one(\baru(k)-x^*)^\T,
\end{array}
\]
it follows that
\begin{equation}\label{eq:last_arrow_proof4_alg3}
\begin{array}{rcl}
\|\bq\|_\Fro^{\lambda, K}
&\leq& \|\cbx\|_\Fro^{\lambda, K} + (\sqrt{n})^{-1}\|(\one-\bv)^\T\bx\|_{\Fro}^{\lambda, K} + \sqrt{n}\|\baru - x^*\|_\Fro^{\lambda, K}\\
&\leq& \left(1+\frac{\sqrt{n}}{\lambda}\left(\sqrt{\frac{L(1+\eta)}{\barmu\eta}+\frac{\hatmu}{\barmu}\beta}\right)\right)\|\cbx\|_\Fro^{\lambda, K}+2\sqrt{n} \|\barx(0) - x^*\|_\Fro\\
&    &+\left((\sqrt{n})^{-1}+\frac{1}{\lambda}\left(\sqrt{\frac{L(1+\eta)}{\barmu\eta}+\frac{\hatmu}{\barmu}\beta}\right)\right)\|(\one-\bv)^\T\bx\|_{\Fro}^{\lambda, K}.
\end{array}
\end{equation}
Finally, we can bound the last term in \eqref{eq:last_arrow_proof4_alg3} as follows
\begin{equation}
\begin{array}{rcl}\label{eq:last_arrow_proof5_alg3}
\|(\one-\bv(k))^\T\bx(k)\|_\Fro
& =  &\|(\one-\bv(k))^\T(I-\frac{1}{n}\one\one^\T)\bx(k)\|_\Fro\\
&\leq&\sqrt{n^2-n}\|\cbx(k)\|_{\Fro}\\
&\leq&n\|\cbx(k)\|_{\Fro}.
\end{array}
\end{equation}
Substituting \eqref{eq:last_arrow_proof5_alg3} into \eqref{eq:last_arrow_proof4_alg3} gives \eqref{eq:last_arrow_proof0_alg3}.
\end{proof}

\subsection{Linear Convergence of Push-DIGing}

Next, we provide a convergence rate estimate for the Push-DIGing algorithm.
\smallskip
\begin{theorem}[\textbf{Algorithm 2: Explicit geometric rate over time-varying directed graphs}]\label{theorem:final_bound_alg3} Suppose Assumptions \ref{ass:smooth}, \ref{ass:strongly_convex}, \ref{ass:connectedness2}, and \ref{ass:matrix_C} hold. Let $B$ be a large enough integer constant such that
\[
\sig\triangleq Q_1\left(1-\frac{1}{n^{(2+n\BG)n\BG}}\right)^{\frac{B-1}{n\BG}}<1.
\]
Also define the constant $J_2$ as follows:
\[
\Cnpush=3Q_1\|\bV^{-1}\|_{\max}^1\barkappa B(\sig+Q_1(B-1))(1+\sqrt{n})\left(1+4\sqrt{n}\sqrt{\barkappa}\right).
\]
Then, for any step-size $\alpha\in\left(0,\frac{1.5(1-\sig)^2}{\barmu\Cnpush}\right]$, the sequence $\{\bx(k)\}$ generated by Push-DIGing converges to the unique optimal solution $\bx^*$ at a global R-linear (geometric) rate $O(\lambda^k)$, where $\lambda$ is given by
\[\label{eq:final_rate_alg3}
\lambda=\left\{
         \begin{array}{ll}
           \sqrt[2B]{1-\frac{\alpha\barmu}{1.5}}, &\text{ if } \alpha\in\left(0,\frac{1.5\left(\sqrt{\Cnpush^2+(1-\sig^2)\Cnpush}-\sig\Cnpush\right)^2}{\barmu\Cnpush(\Cnpush+1)^2}\right],\\
           \sqrt[B]{\sqrt{\frac{\alpha\barmu\Cnpush}{1.5}}+\sig}, &\text{ if } \alpha\in\left(\frac{1.5\left(\sqrt{\Cnpush^2+(1-\sig^2)\Cnpush}-\sig\Cnpush\right)^2}{\barmu\Cnpush(\Cnpush+1)^2},\frac{1.5(1-\sig)^2}{\barmu\Cnpush}\right].
         \end{array}
       \right.
\]
\end{theorem}
\smallskip
\begin{proof}
The proof is similar to that of Theorem \ref{theorem:final_bound}. Specifically, we collect all the gains as follows:
\[
\begin{array}{ll}
\text{i)}\ \gamma_1=L\left(1+\frac{1}{\lambda}\right);
&\quad\text{ii)}\
\gamma_2=
\frac{Q_1\|\bV^{-1}\|_{\max}^1\lambda(1-\lambda^B)}{(\lambda^B-\sig)(1-\lambda)};\\
\text{iii)}\ \gamma_3=\frac{\alpha}{\lambda^B-\sig}\left(\sig+Q_1\frac{1-\lambda^{B-1}}{1-\lambda}\right);
&\quad\text{iv)}\ \gamma_4=(1+\sqrt{n})\left(1 + \frac{\sqrt{n}}{\lambda}\sqrt{\frac{L(1+\eta)}{\barmu\eta}+\frac{\hatmu}{\barmu}\beta}\right).
\end{array}
\]
To apply the small gain theorem (Theorem \ref{theorem:small}), we need $\gamma_1\gamma_2\gamma_3\gamma_4<1$, that is,
\[\label{eq:final_1_alg3}
\begin{array}{l}
\frac{Q_1\|\bV^{-1}\|_{\max}^{1}\alpha L(1+\lambda)}{(\lambda^B-\sig)^2}\left(\frac{1-\lambda^B}{1-\lambda}\right)\left(\sig+Q_1\frac{1-\lambda^{B-1}}{1-\lambda}\right)(1+\sqrt{n})\left(1 + \frac{\sqrt{n}}{\lambda}\sqrt{\frac{L(1+\eta)}{\barmu\eta}+\frac{\hatmu}{\barmu}\beta}\right)<1.
\end{array}
\]
The other restrictions on $\alpha$, $\beta$, $\eta$ and $\lambda$ are the same as those in \eqref{eq:final_2}, \eqref{eq:final_3}, \eqref{eq:final_4}, and \eqref{eq:final_5}. Choosing $\beta=2L/\hatmu$ and $\eta=1$, and further using $0.5\leq\lambda<1$ and $(1-\lambda^{B-1})/(1-\lambda)\leq B-1$, relations \eqref{eq:final_1_alg3} and \eqref{eq:final_5} together imply that
\begin{equation}\label{eq:final_proof1_alg3}
\alpha\leq\frac{(\lambda^B-\sig)^2}{2Q_1\|\bV^{-1}\|_{\max}^1 LB(\sig+Q_1(B-1))(1+\sqrt{n})\left(1+4\sqrt{n}\sqrt{\barkappa}\right)}.
\end{equation}

Next, similar to the proof of Lemma \ref{lemma:bound_R-linear} [cf.\ \eqref{eq:final_proof5}], it remains to show that there exists some $\lambda\in(\sqrt[B]{\sig},1)$ such that
\begin{equation}\label{eq:final_proof3_alg3}
\left[\frac{1.5(1-\lambda^{2B})}{\barmu},\frac{1.5(\lambda^B-\sig)^2}{\barmu\Cnpush}\right]\neq\emptyset,
\end{equation}
where $\Cnpush=3Q_1\|\bV^{-1}\|_{\max}^1\barkappa B(\sig+Q_1(B-1))(1+\sqrt{n})\left(1+4\sqrt{n}\sqrt{\barkappa}\right)$. Noticing the similarity between \eqref{eq:final_proof3_alg3} and \eqref{eq:final_proof5}, by the same argument as in the proof of Theorem \ref{theorem:final_bound} (starting from \eqref{eq:final_proof5}), we can obtain the statement of the theorem.
\end{proof}
\smallskip

\section{Numerical Experiments}\label{sec:num_exp}

Consider a decentralized estimation problem: each agent $i\in\{1,\ldots,n\}$ has its own observation $y_i$ given by $y_{i}=M_{i} x+e_{i}$, where $y_{i}\in\R^{m_i}$ and $M_{i}\in\R^{m_i\times p}$ are known data, $x\in\mathbb{R}^p$ is unknown, and $e_{i}\in\mathbb{R}^{m_i}$ is some noise. The goal is to estimate $x$. In this experiment we use the Huber loss, which is known to be robust to outliers, and it allows us to observe both sublinear and linear convergence. The corresponding optimization problem is: $\min_{x\in\R^p} f(x)=\frac{1}{n}\sum_{i=1}^n\left\{\sum_{j=1}^{m_i}H_\xi(M_{i,j} x-y_{i,j})\right\}$, where $M_{i,j}$ and $y_{i,j}$ are the $j$-th row of matrix $M_{i}$ and vector $y_{i}$, respectively. The Huber loss function $H_\xi$ is defined by
\[
H_\xi(a)=\left\{
           \begin{array}{cl}
             \frac{1}{2}a^2, &\text{ if $|a|\leq\xi$ ($\ell_1$ zone)},\\
             \xi(|a|-\frac{1}{2}\xi), &\text{ otherwise ($\ell_2^2$ zone)}.
           \end{array}
         \right.
\]
In all experiments, we set $n=12$, $m_i=1$ for all $i$, and $p=3$. The data $M_{i}$, as well as the noise $e_{i}$, $\forall i$, are firstly generated following the standard normal distribution and then re-normalized so that $L_i=1$ for all $i$. The parameter $\xi$ is set to $2$. The algorithm starts from $x_{i}(0)=0$ for all $i$. We scale $e_{i}$ so that $x_{i}(0)$ is located in the $\ell_1$ zone for all agents $i$. The optimal solution $x^*$ is randomly but artificially selected so that $x^*$ is in the $\ell_2^2$ zone and $\|x^*-x_i(0)\|=300,\ \forall i$. \mr{Note that while the objective is restricted strongly convex, not all individual functions are strongly convex.}

\begin{figure}[h]
\begin{center}
\subfigure{
\begin{minipage}[h]{0.45\textwidth}
\includegraphics[height=9em]{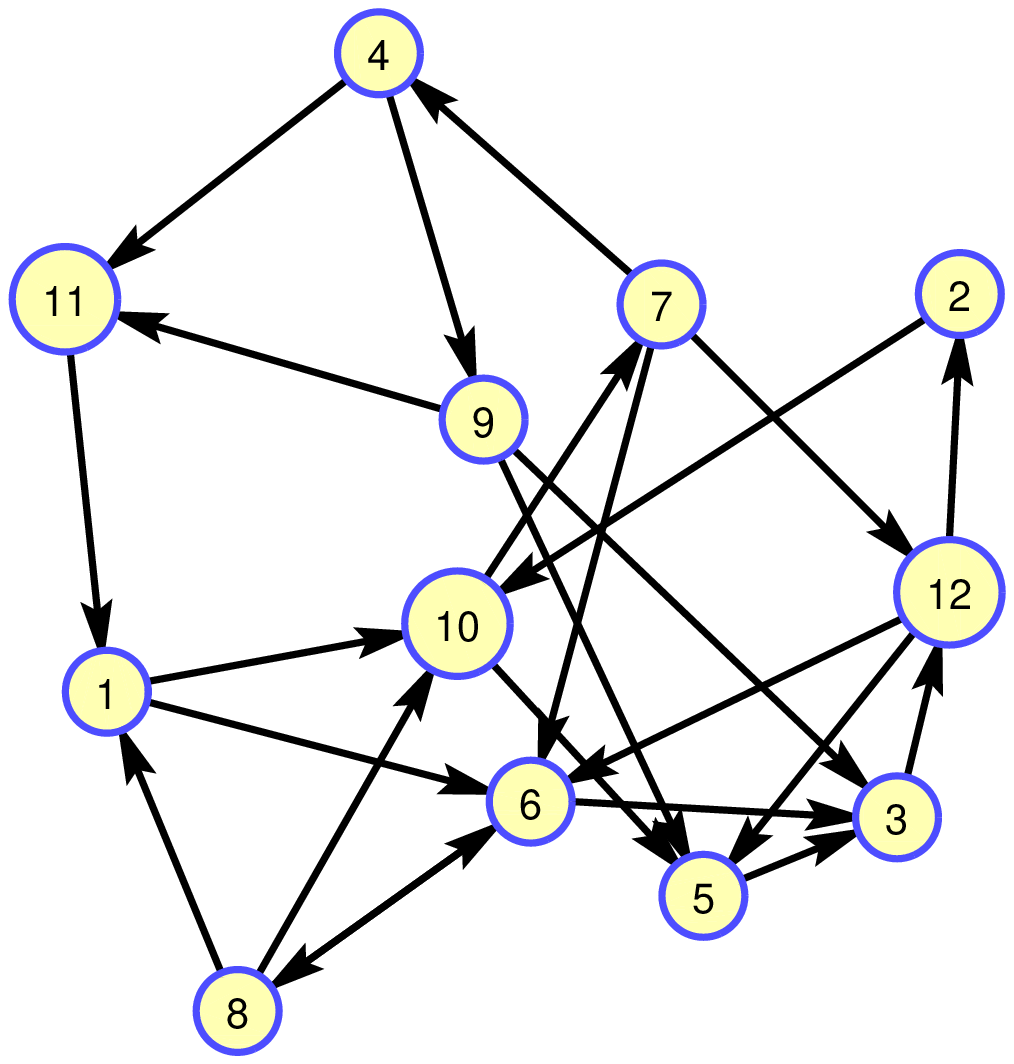}\\
\includegraphics[height=9em]{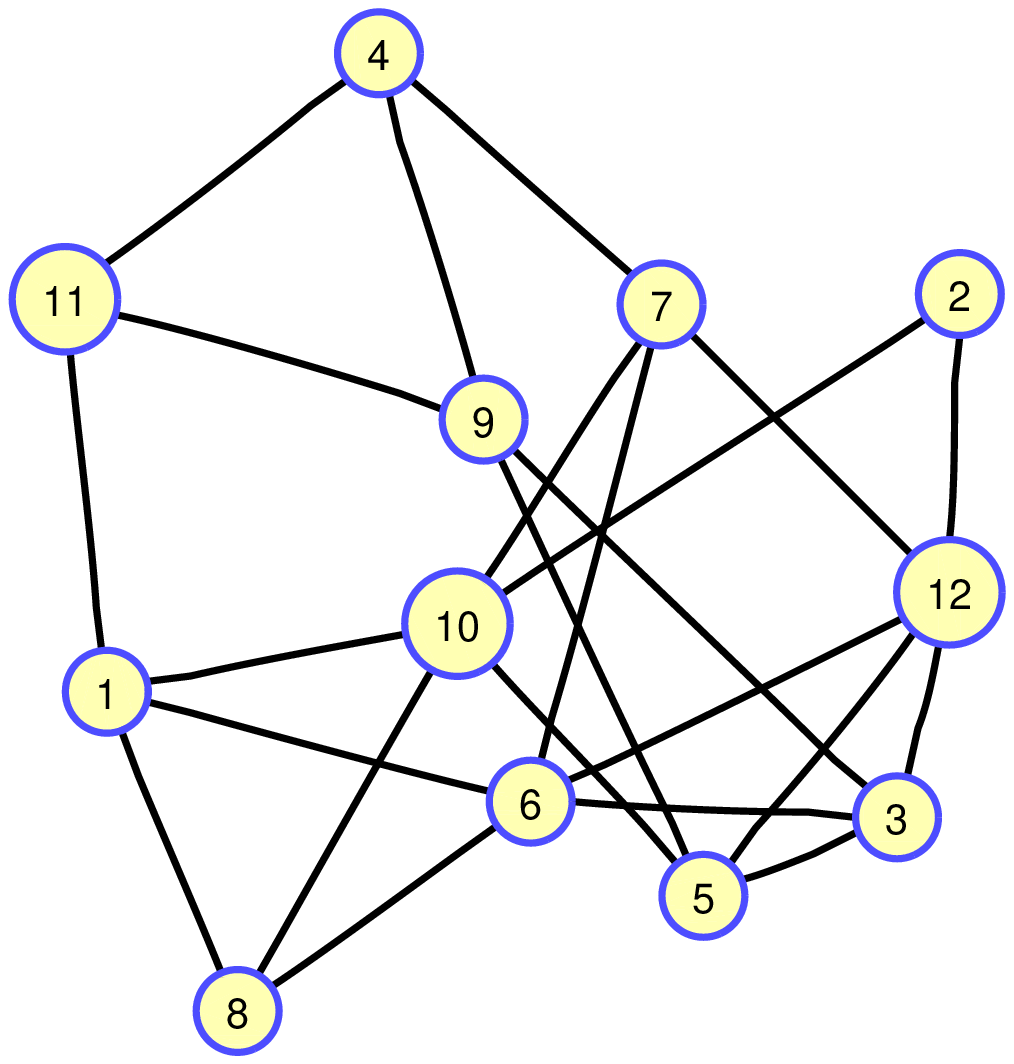}
\end{minipage}
}\hspace{-7em}
\subfigure{
\begin{minipage}[h]{0.55\textwidth}
\includegraphics[height=18em]{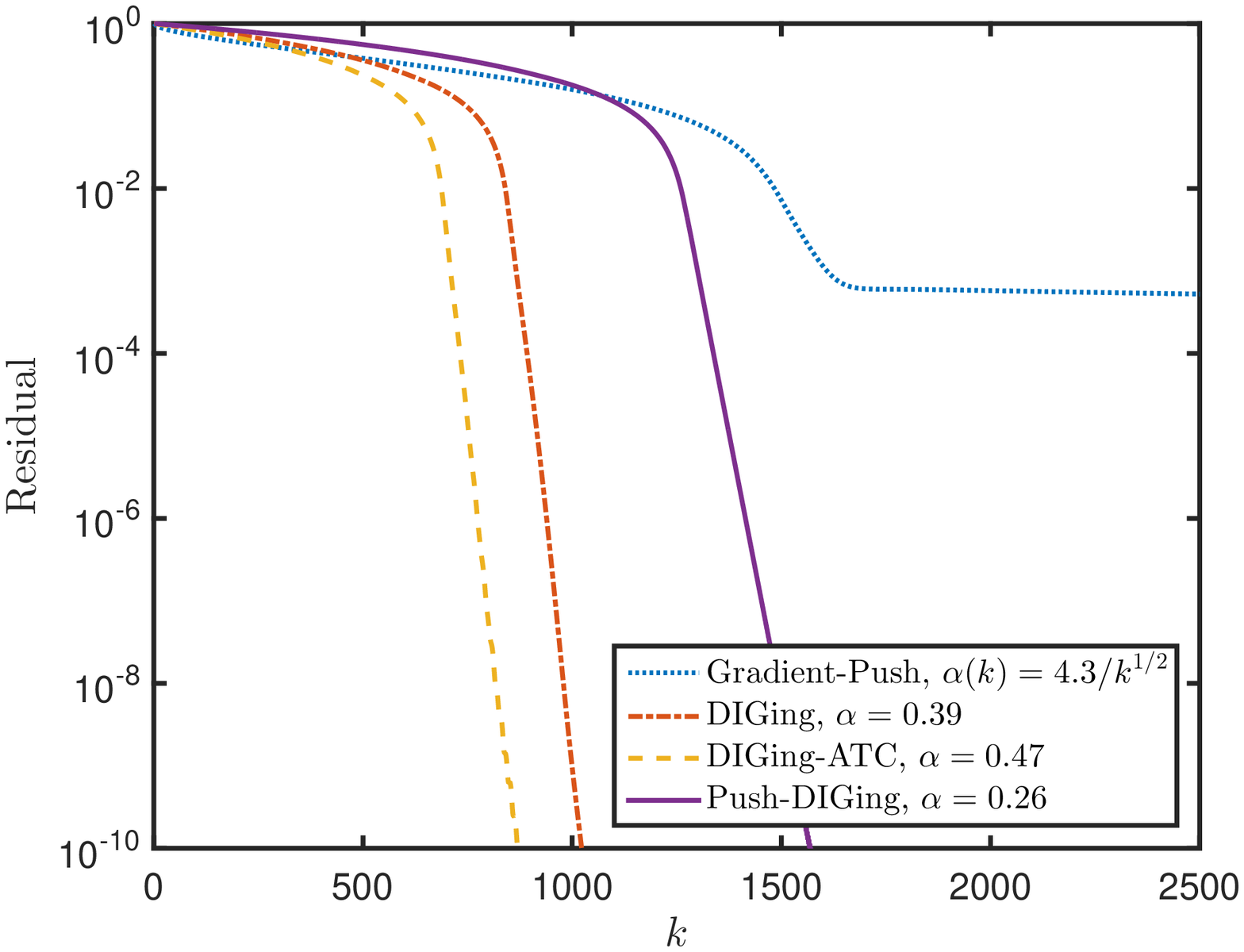}
\end{minipage}
}\caption{The plots are showing the underlying directed and undirected graphs for experiments. The plot to the right shows the residuals $\frac{\|\bx(k)-\bx^*\|_\Fro}{\|\bx(0)-\bx^*\|_\Fro}$
for the time-invariant directed graph illustrated on the top-left. Step-sizes have been \emph{hand-optimized} to give faster convergence and more accurate solution for all algorithms.}\label{eps:RLS_TI}
\end{center}
\end{figure}
\begin{figure}[h]
\begin{center}
\includegraphics[height=15em]{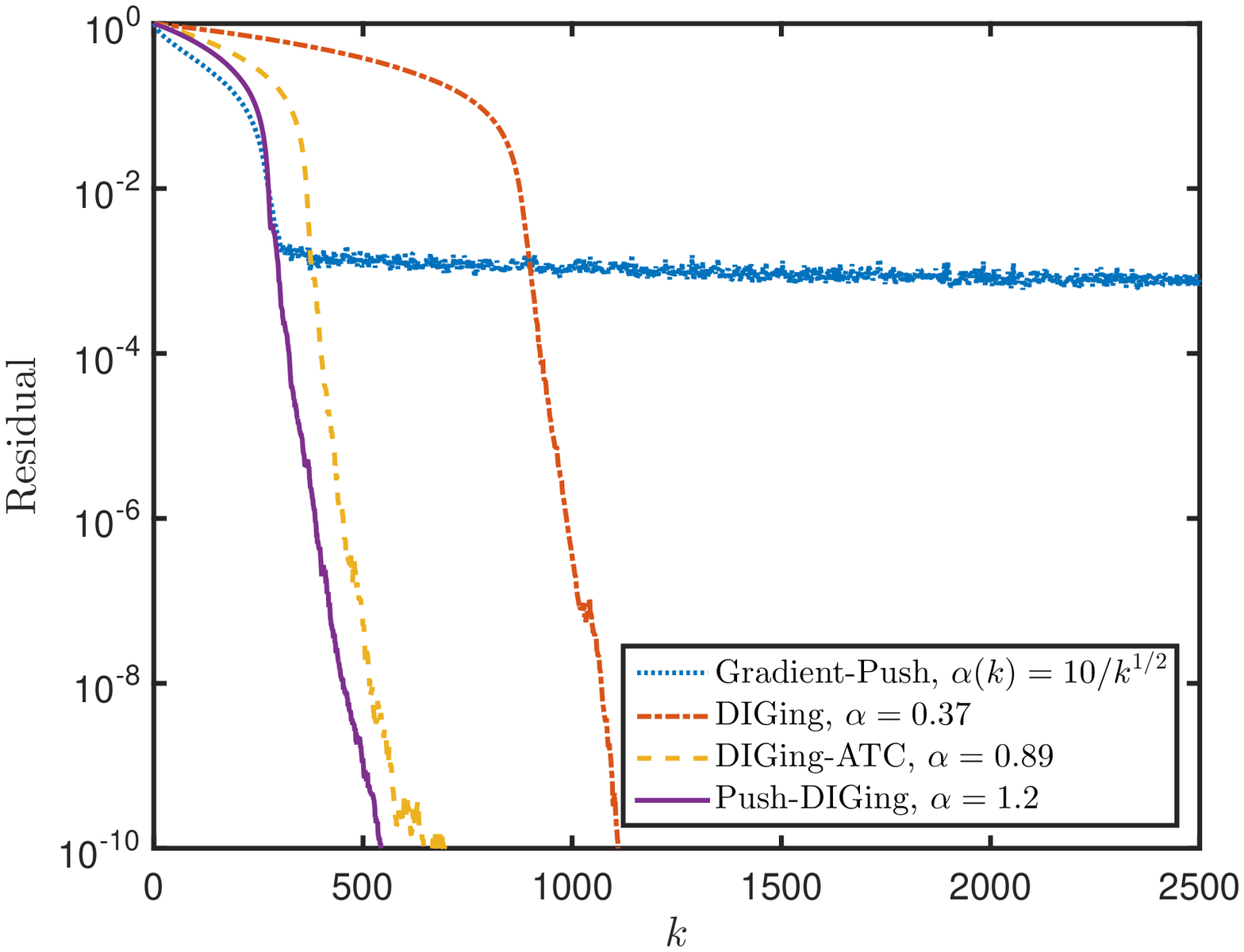}\hspace{0.2em}
\includegraphics[height=15em]{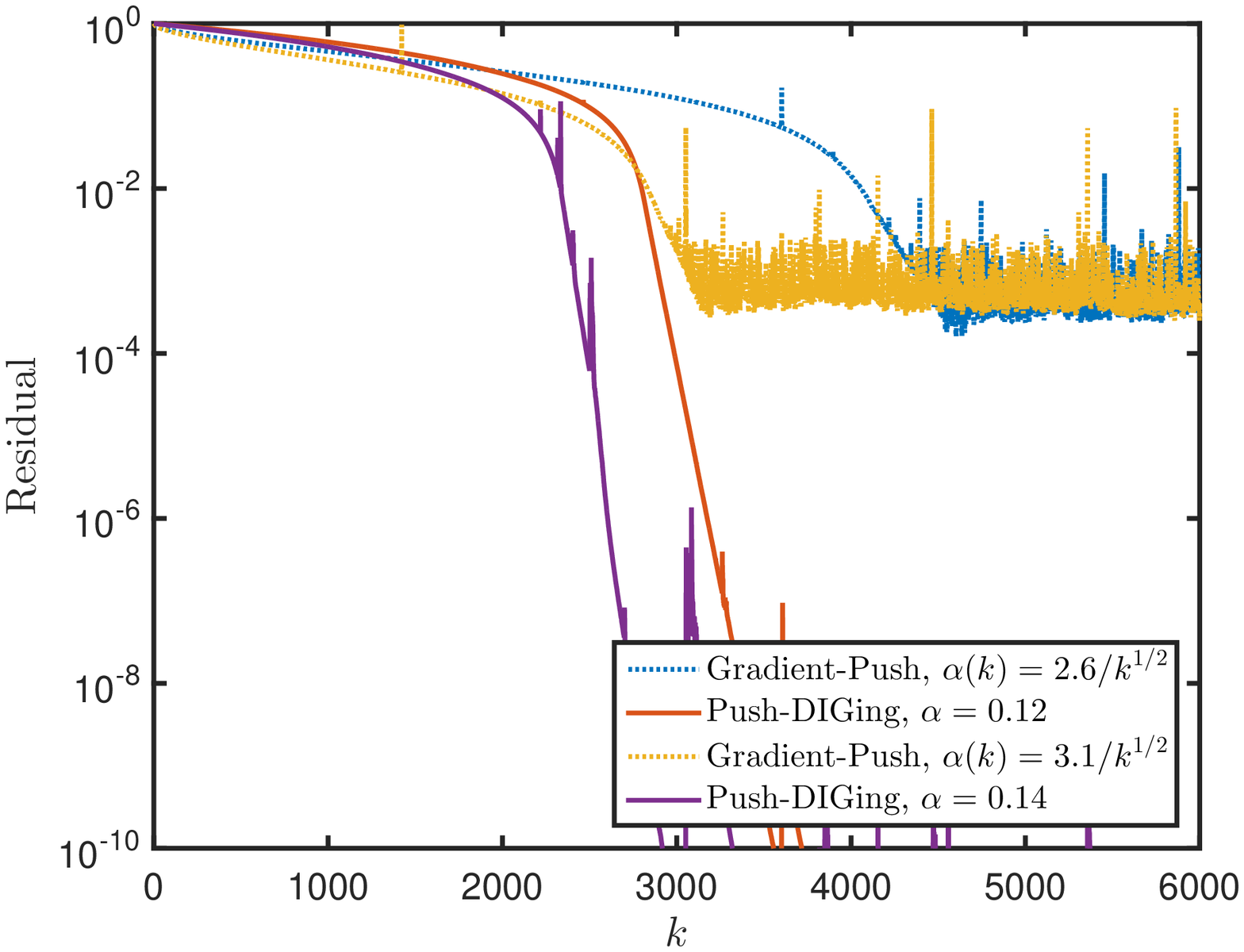}
\caption{The plots to the left and to the right are showing the residuals $\frac{\|\bx(k)-\bx^*\|_\Fro}{\|\bx(0)-\bx^*\|_\Fro}$ for a time-varying undirected graph sequence and a time-varying directed graph sequence, respectively. Step-sizes have been \emph{hand-optimized} to give faster convergence and more accurate solution for all algorithms.}\label{eps:RLS_TV}
\end{center}
\end{figure}

We conduct three experiments on (i) Time-invariant directed graphs; (ii) Time-varying undirected graphs; and (iii) Time-varying directed graphs. For (i), the underlying directed graph $\GTIdir=\{\V,\A\}$ is illustrated in the top left of Fig. \ref{eps:RLS_TI}. $\GTIdir$ is randomly generated with $24$ arcs and $12$ vertices, and guaranteed to be strongly connected. The mixing matrix $\Wdo$ is chosen as the fastest (asymmetric) linear consensus matrix \cite{Xiao2004}. For (ii), the underlying undirected graph $\GTIun=\{\V,\E\}$ is generated by simply taking out the direction of every arc of $\GTIdir$ (see the bottom left of Fig. \ref{eps:RLS_TI}). $\GTIun$ has $23$ edges in total. At iteration $k$, we generate the graph instance of the time-varying undirected graph $\GTVun(k)=\{\V,\E(k)\}$ through randomly uniformly sampling $\E(k)$ from $\E$ with $40$ percent. The mixing matrix $\Wdo(k)$ is generated through Metropolis weights. For (iii), we generate the graph instance of the time-varying directed graph $\GTVdir(k)=\{\V,\A(k)\}$ through randomly uniformly sampling $\A(k)$ from $\A$ of $\GTIdir$ with $80$ percent. In all experiments, the mixing matrix $\Wco(k)$ for Push-DIGing is generated based on the out-degree information \cite{Nedic2015}.

We plot push-gradient method \cite{Nedic2014}, DIGing, DIGing-ATC (Aug-DGM with step-size matrix $\alpha I$), and Push-DIGing in Fig. \ref{eps:RLS_TI} and Fig. \ref{eps:RLS_TV}. In the experiments, we observe that DIGing and its variants all have R-linear rates while push-gradient method only has sublinear rate even if the objective is smooth and strongly convex. \mr{Since our analyses in the above theorems and corollaries have all been worst-case analyses, the bounds on step-sizes we have given are often conservative. We therefore choose to use \emph{hand-optimized} step-sizes for the numerical experiments.}

\mr{We first discuss simulation performed over an undirected communication graph. In this case, we have also tested the EXTRA algorithm from reference \cite{Shi2015}. Interestingly, we find the convergence curve of the EXTRA algorithm almost identical\footnote{Since the two curves are almost identical, EXTRA is not plotted in the figure.} to that of the DIGing algorithm when the same step size is chosen for both algorithms. As commented in Remark \ref{remark:DIGing-ATC}, we have also conducted the numerical experiments for DIGing-ATC and indeed find it faster, in terms of  number of iterations, than DIGing (see the left sub-figure of Fig. \ref{eps:RLS_TV}). However, note that the per-iteration communication cost of DIGing-ATC is higher.}

\mr{For an undirected graph, Push-DIGing reduces to DIGing with partial ATC strategy which only needs one round of communication per iteration. We also plot Push-DIGing for undirected graphs in the left of Fig. \ref{eps:RLS_TV} and find that Push-DIGing can be as fast as DIGing-ATC in terms of number of iteration but at actually at a half communication cost per iteration. We speculate that this might be because after using the ATC strategy once (in Push-DIGing), the bound of step size/convergence rate has reached the limit of centralized gradient descent and no more improvement can be obtained with  ATC structure.}

\mr{Although DIGing and DIGing-ATC are designed for undirected graphs, we also test them over directed graphs using an asymmetric doubly stochastic matrix. The results are shown in the right sub-figure of Fig. \ref{eps:RLS_TI}.  Note that constructing an asymmetric doubly stochastic matrix over a directed graph will require a graph balancing algorithm with considerable overhead. We note that the convergence curve of the DEXTRA/ExtraPush algorithm from reference \cite{Xi2015,Zeng2015} is not plotted in Fig. \ref{eps:RLS_TI} since we did not find a stable step size for the algorithm in this case. This might be because the individual functions are not strongly convex. }
	
\mr{Finally, for time-varying directed graphs, only Push-DIGing is plotted in the right sub-figure of Fig. \ref{eps:RLS_TV} as this is the only algorithm over time-varying directed graphs with geometric convergence. We do not plot DIGing/DIGing-ATC since they need real time graph balancing which is not practically implementable when the graph varies.}

\section{Conclusion}\label{sec:concl}

In this paper, we considered a class of protocols for distributed optimization based on the idea of ``distributed inexact gradient'' and ``gradient tracking''. Under strong convexity, we studied the convergence rates of the algorithms over time-varying directed/undirected graphs. Using the small gain theorem, we showed that our protocols converge at some global R-linear rates for strongly convex functions, and were able to obtain explicit bounds on the rates.

An open question is to obtain improved estimates on the convergence rates of our method, especially as far as scaling with the number of agents, $n$, goes. Furthermore, extensions to more complex optimization models containing local constraints, couplings among agents in the objectives would be of considerable interest.

\section*{Acknowledgement}

We thank C\'{e}sar A. Uribe and Thinh T. Doan for their helpful discussions.
\bibliographystyle{siamplain}
\bibliography{document}
\end{document}